\documentclass{article}
\usepackage{amssymb}
\usepackage{amsfonts}
\usepackage{amsmath}
\usepackage{color}

\newtheorem{theorem}{Theorem}[section]

\newtheorem{definition}[theorem]{Definition}

\newtheorem{lemma}[theorem]{Lemma}

\def\ud{\, \mathrm{d}}
\numberwithin{equation}{section}

\begin{document}

\title{A Bond Option Pricing Formula in the Extended CIR Model, with an
Application to Stochastic Volatility}
\author{Zheng Liu\footnote{Morgan Stanley, New York city, E-mail: Zheng.Liu@cgu.edu}, Qidi Peng\footnote{Institute of Mathematical Sciences, Claremont Graduate University, E-mail: Qidi.Peng@cgu.edu} and Henry Schellhorn\footnote{Corresponding author, Institute of Mathematical Sciences, Claremont Graduate university, E-mail: Henry.Schellhorn@cgu.edu}}
\maketitle

\begin{abstract}
We provide a complete representation of the interest rate in the extended
CIR\ model. Since it was proved in Maghsoodi (1996) that the representation
of the CIR process as a sum of squares of independent Ornstein-Uhlenbeck
processes is possible only when the \textit{dimension} of the interest rate
process is integer. We use a slightly different representation, valid when
the dimension is not integer. Our representation uses an infinite sum of
squares of basic processes. Each basic process can be described as an
Ornstein-Uhlenbeck process with jumps at fixed times. In this case, the
price of a bond option resembles the Black-Scholes formula, where the normal
distribution is replaced by a new distribution, which generalizes the
non-central chi-square distribution. For practical purposes, the bond price
is however better calculated by inverting a Laplace transform. We generalize
the result to a model where volatility is stochastic.
\end{abstract}

\section{Introduction}

The Cox-Ingersoll-Ross term structure model, in short the CIR model, was
introduced in Cox et al. (1985a), (1985b). In this model, the spot interest
rate $r(t)$ is assumed to follow a squared Bessel process:%
\begin{eqnarray}
\label{CIR}
\ud r(t) &=&(-b(t)r(t)+\theta (t))\ud t+\sigma (t)\sqrt{r(t)}\ud B(t)
\label{masterstoch2} \\
r(0) &=&r_{0}>0,  \notag
\end{eqnarray}

where $B$ is a standard Brownian motion. In the original specification of Cox et al. (1985a), (1985b), the speed of mean reversion $b\geq 0$, the volatility $%
\sigma >0$ and the parameter $\theta >0$ are assumed constant. In this
paper, we consider the case where all these parameters are continuous
functions of time. This is sometimes referred to as the \textit{extended
CIR\ model}.

Several features of the CIR\ model are particularly attractive. Firstly, it
can be justified by general equilibrium considerations, see Cox et al.
(1985a). Secondly, the interest rate is always positive and stationary. Cox
et al. found that its distribution follows a non-central chi-square
distribution. Finally, there is a closed form formula for the bond price.
For practitioners however, the main shortcoming of the constant parameters
version of the model is that it cannot reproduce the original term structure
of interest rates. This fact was highlighted by several authors (Hull
(1990), Keller-Ressel and Steiner (2008), Yang (2006) and all the references
therein): yield curves can be only normal, inverse, or humped. The extended
CIR model, however, has enough parameters to be fitted to the original yield
curve.

Maghsoodi (1996), Jamshidian (1995), and Rogers (1995) propose a
representation of the extended CIR\ model as a sum of squares of
Ornstein-Uhlenbeck processes when the dimension $d\equiv 4\theta (t)/\sigma
^{2}(t)$ is constant and integer. As a consequence, the interest rate
follows the generalized chi-square distribution. Maghsoodi (1996), and
Shirakawa (2002) also propose a representation of the interest rate as a
time-changed lognormal process. However, as the latter author states, "it is
difficult to derive the probability distribution of the squared Bessel
processes with time-varying dimensions explicitly". Brigo and Mercurio
(2006) state that no solution to that problem has been found. We view this
as the principal contribution of our paper. While the characteristic
function of the interest rate has a remarkably simple structure, there are
many different ways to invert it to characterize the probability
distribution. One of them is by moment approximation. In the constant
parameters case, the bond price is determined by solving a Riccati equation.
There is in general no closed form formula for the bond price, and authors
generally claim that their solution is in closed form "up to solution of an
ODE". We note that Mannoli et al. (2008) also study the extended
CIR\ model, and provide a closed form solution for the bond price as an
infinite series. We then generalize Maghsoodi's result (1996) for the price
of a bond option to the case where the dimension $d(t)$ is any positive
continuous function.

There are alternate approaches to extend the CIR model. Brigo and Mercurio
(2001) find that a deterministic shift of the CIR model is analytically
tractable. The obvious drawback of making the parameters of the CIR model
functions of time is the problem of overparameterization: the parameters
will not be robust in a change of regime. Several authors consider instead a
CIR\ model of interest rates with constant parameters and stochastic
volatility. For instance, Longstaff and Schwartz (1992), and Duffie and Kan
(1996), consider a generalized two-factor CIR\ model, where one factor is
the interest rate, and the other one is its volatility. The volatility in
that model is a variation of the volatility in the popular Heston (1993)
model. Cotton, Fouque, Papanicolaou and Sircar (2004) calculate an
asymptotic expansion of the bond price in such a model (with constant
parameters) when the speed of mean reversion is fast. Fouque and Lorig
(2011) generalize this model to a model with a volatility of volatility.

Finally, we note that several authors have generalized the CIR model in a
different way using multiple factors. We refer the reader to the references
contained in Chen, Filipovic, and Poor (2004) and Gourieroux and Monfort
(2011).

We then incorporate stochastic volatility $v(t)\ $in our extended CIR model.
If $b$ and $\theta $ are deterministic, and $v(t)$ is independent of $r(t)$
the distribution of the interest rate, conditional on volatility, is the
same as in the deterministic volatility case. As an illustration of our
approach, we give a semi-closed formula of the price of a bond option when
the price of the underlying bond at time $t$ can be expressed
a function as a function of $r(t),v(t)$ and time $t$. We also consider a
model where volatility is correlated with the interest rate.

In the same way that it is not too difficult to generalize the one-factor
CIR model to multiple factors (see Duffie and Kan (1996)), we believe it is
not difficult to generalize our results on the extended CIR model to
multiple factors. However, we leave this for future research.

The structure of this paper is as follows. In Section 2, we consider the
case where volatility is deterministic, namely the extended CIR model
\textit{per se}. In Section 3, we provide an analytic formula (up to
solution of ODEs) for a bond option with deterministic volatility. In
Section 4, we extend our results to stochastic volatility.

\section{The Extended CIR\ Model}
\begin{definition}
\label{def1}
We let $T>0$ be the maturity of the bond underlying
the bond option.
\end{definition}

\begin{definition}
\label{def2}
We define the number of degrees of freedom (also
called \textit{dimension} in Maghsoodi)
\begin{equation*}
d(t)=\frac{4\theta (t)}{\sigma ^{2}(t)}.
\end{equation*}%
\end{definition}

Since $\theta $ is assumed to be strictly positive, and $\sigma $ non-zero, $%
d(t)$ is also strictly positive.

\bigskip

\subsection{An Informal View of Our Representation}

\subsubsection{The Rate as a Sum of Squared Gaussian Processes: an
Incomplete, but Useful Initial Concept}

Our goal is to generalize the representation of the rate $r(t)$ as a sum with a
constant number of terms of squares of Ornstein-Uhlenbeck processes into a
sum of time-varying numbers of such terms. Heuristically, it is easier to
represent this sum as an integral. Before introducing the rigorous
definitions, we explain informally the method. We introduce a
collection $\left\{ W(.,u)\right\} $ of independent Brownian motions
\footnote{%
Rather than starting with a countably infinite collection of processes $%
x(.,u)$ we will start with a countable finite collection of processes in the
formal proof $\{\tilde{x}_{i}(.)\}$. See below.}. We write the Gaussian
processes as:%
$$
\frac{\partial x(t,u)}{\partial t}=\frac{1}{2}(-b(t)+\text{"plug"})x(t,u)\ud t+\frac{\sigma (t)}{2}%
\ud W(t,u)  \label{dtx}
$$
for some unknown functional "plug". In particular, if we take the representation:%
$$
r(t)=\int_{0}^{d(t)}x(t,u)^{2}\ud u.  \label{reprofrate}
$$

Then informally, taking derivatives results in:%
$$  \label{dr}
\ud r(t)=\int_{0}^{d(t)}\frac{\partial (x(t,u)^{2})}{\partial t}\ud u+x(t,d(t))^{2}d'(t).
$$

We see that:%
\begin{eqnarray*}
\ud r(t) &=&(-b(t)r(t)+\frac{\sigma ^{2}(t)d(t)}{4})\ud t+\sigma (t)\sqrt{r(t)}%
\ud B(t)+\int_{u=0}^{d(t)}x(t,u)^{2}\text{"plug"}\ud u \\
&&+x(t,d(t))^{2}d'(t).
\end{eqnarray*}

If the last two terms cancel out, by Definition \ref{def1}, we obtain the right
equation for $r$. Intuitively, the last two terms will cancel out if the
"plug" term includes a Dirac delta function:%
\begin{equation*}
\text{"plug"}=-\delta (u-d^{-1}(t))d'(t),
\end{equation*}
provided $d$ is invertible.

\subsubsection{Application of the Infinite Divisibility of the Chi-Squared
Distribution}

We want to represent $\int_{0}^{d(t)}x(t,u)^{2}\ud u$ in (\ref{reprofrate}) as
a sum. To this end we discretize the $x(t,u)$ process in the $u$ direction and replace $%
\int_{0}^{d(t)}x(t,u)^{2}\ud u$ by a sum of squares of basic processes $\{\tilde{x%
}_{i}(t)\}_i$. Since the "number" of terms $d(t)$ is real-valued, an idea is to
approximate $d(t)$ by the ratio $\frac{[nd(t)]}{n}$. We then construct a sum
of $[nd(t)]$ terms, which we call \footnote{%
The subscript $_{low}$ indicates that we truncate the integral into a lower
value.} $Z_{low}^{(n)}(t)$:

\begin{equation*}
Z_{low}^{(n)}(t)=\sum_{i=1}^{[nd(t)]}\tilde{x}_{i}(t)^{2}.
\end{equation*}%
We should then "divide" $Z_{low}^{(n)}(t)$ by $n$, and then make $n$ go to
infinity. The reason why we put quotes around "divide" is that the variance
of $Z_{low}^{(n)}(t)/n$ tends to zero, and thus $Z_{low}^{(n)}(t)/n$ is not
the right representation for the rate. Rather, if $Z_{low}^{(n)}(t)$ had a
non-central chi-squared distribution, we could represent $Z_{low}^{(n)}$ as
a sum of i.i.d. processes $\{r_{j,low}^{(n)}(t)\}_j$
$$
Z_{low}^{(n)}(t)=\sum_{j=1}^{n}r_{j,low}^{(n)}(t).  \label{Zlow}
$$

\subsubsection{The Rate as a Double Limit}

It turns out that $Z_{low}^{(n)}(t)$ does not have the scaled non-central
chi-squared distribution. However, it is possible to define a process $%
r_{j,mid}^{(n,M)}(t)$ so that, conditionally on the information available at
time $\frac{(m-1)}{M}T$, the random variable $r_{j,mid}^{(n,M)}(\frac{m}{M}%
T)\ $has approximately the scaled non-central chi-squared distribution.
Since the latter is computable, we can compute the characteristic function
of $r_{j,mid}^{(n,M)}(\frac{m}{M}T)$ conditionally on $r_{j,mid}^{(n,M)}(%
\frac{(m-1)}{M}T)$. We then iterate this backward induction by calculating
this characteristic function conditionally on $r_{j,mid}^{(n,M)}(\frac{(m-2)%
}{M}T)$, and so far until time zero. We will then take a double limit of
this characteristic function, by making both $n$ and $M$ go to infinity.

\subsection{Definitions}

Let $(\Omega ,\mathcal{F}_t,\mathbb{P})$ be a probability space where the
filtration $\mathcal{F}_{t}$ is generated by a countably infinite collection
$\{W_{i}\}$ of independent Brownian motions. The probability measure $%
\mathbb{P}$ is the risk-neutral measure.

\bigskip

\textbf{Observation}

Most subsequent definitions will not be needed to understand part (i) of
Theorem \ref{thm2}, and the reader may skip this tedious part, which will be
necessary only in the stochastic volatility section.

\bigskip

We start by defining

\begin{eqnarray}
&&d^{(n)}(t) =nd(t);  \notag \\
&&M_{d}^{(n)} =\max_{t\in \lbrack 0,T]}d^{(n)}(t);  \notag \\
&&t_{m}^{M} =\frac{m}{M}T\text{ \ \ \ \ }1\leq m\leq M.  \label{t_m_M}
\end{eqnarray}

Definition (\ref{t_m_M}) will be a standing definition, in the sense that,
for each $m,M$, (\ref{t_m_M}) will hold throughout the paper.

\bigskip

\textbf{Assumption 1}

We suppose that, over the interval $[0,T]$, $d(t)$ satisfies Dirichlet's
condition, i.e., it is differentiable, the number $Q$ of minima and maxima of $%
d(t)$ is finite, and that the set of critical points of $d$ (i.e. the points
where $d'(t)=0$) has measure zero. As mentioned above, $d(t)$ is strictly
positive. Also, we assume that $\theta(t)$, $b(t)$ and $\sigma(t)$ are positive real-valued continuous
functions on $[0,T]$. This ensures that (\ref{CIR}) has a pathwise unique strong solution
(see Magshoofi (1996)).

\subsubsection{Definition of $\tilde{x}_{i}^{(n)}(t)$}

We have to consider separately each decreasing and increasing branch of $%
d^{(n)}(t)$. We now have 3 subcases. Accordingly, the processes $\tilde{x}%
_{i}^{(n)}(t)$ can at any time $t$ either:

\begin{itemize}
\item not jump;

\item jump up by a half;

\item jump down by a half.
\end{itemize}

Let $\mathcal{J}_{down}^{(n)}(i)$ for $i\geq 1$ be the set of solutions of:%
\begin{equation*}
d^{(n)}(t)=i\text{ \ such that}\frac{\ud d^{(n)}(t)}{\ud t}>0,\text{ \ \ }t\in \lbrack
0,T].
\end{equation*}%
Let $\mathcal{J}_{up}^{(n)}(i)$ for $i\geq 1$\ be the set of solutions of:%
\begin{equation*}
d^{(n)}(t)=i\text{  such that }\frac{\ud d^{(n)}(t)}{\ud t}<0,\text{ \ \ }t\in
\lbrack 0,T].
\end{equation*}

We do not consider the case $\frac{\ud d^{(n)}(t)}{\ud t}\big|_{t}=0$ and $d^{(n)}(t)=i$
\footnote{%
It is fairly easy to note that the cases where $t$ is such that $\frac{%
\ud d^{(n)}(t)}{\ud t}\big|_{t}=0$ and $d^{(n)}(t)=i$ play no role in
the proof (provided that the set of points where $d$ is flat has measure
zero). From now on, we continue the development without taking care of that case,
in order to save space.}. Denote by $N_{up}^{(n)}(i)$ ($N_{down}^{(n)}(i)$)
the cardinality of the set $\mathcal{J}_{up}^{(n)}(i)$ ($\mathcal{J}%
_{down}^{(n)}(i)$), and designate the elements of these sets followingly:%
\begin{eqnarray*}
\mathcal{J}_{up}^{(n)}(i) &=&\{J_{k,up}^{(n)}(i)|0\leq k\leq
N_{up}^{(n)}(i)\}, \\
\mathcal{J}_{down}^{(n)}(i) &=&\{J_{p,up}^{(n)}(i)|0\leq p\leq
N_{down}^{(n)}(i)\},
\end{eqnarray*}

where for convenience we set $J_{0,up}^{(n)}(i)=J_{0,down}^{(n)}(i)=-\infty $%
. Finally, we call the sets of all jumps together with the terminal time:%
\begin{equation*}
\mathcal{J}^{(n)}=\bigcup\limits_{i=1}^{N_{up}^{(n)}(i)}\mathcal{J}%
_{up}^{(n)}(i)\bigcup\limits_{i=1}^{N_{down}^{(n)}(i)}\mathcal{J}%
_{down}^{(n)}(i)\bigcup \{T\}
\end{equation*}

We call $Z_{k,up}^{(n)}(i)$ ($Z_{p,down}^{(n)}(i)$) the first minimizer
(maximizer) of $d^{(n)}$ after $J_{k,up}^{(n)}(i)$ ($J_{p,down}^{(n)}(i)$).
We have then (again, barring the case $\frac{\ud d^{(n)}(t)}{\ud t}=0$ and $%
d^{(n)}(t)=i$ ):%
\begin{equation*}
J_{k,up}^{(n)}(i)<Z_{k,up}^{(n)}(i)<J_{p(k),down}^{(n)}(i)<Z_{p(k),down}^{(n)}(i)<J_{k+1,up}^{(n)}(i),
\label{kandp}
\end{equation*}

where:%
\begin{eqnarray*}
k(p) &=&\left\{
\begin{array}{c}
p \\
p+1
\end{array}%
\right. \text{if}%
\begin{array}{c}
d^{(n)}\text{ is decreasing at }t=0 \\
d^{(n)}\text{ is increasing at }t=0,
\end{array}
\label{kandp_lespetits} \\
p(k) &=&\left\{
\begin{array}{c}
p \\
k-1%
\end{array}%
\right. \text{if}%
\begin{array}{c}
d^{(n)}\text{ is decreasing at }t=0 \\
d^{(n)}\text{ is increasing at }t=0.
\end{array}%
\end{eqnarray*}

For definiteness, we set $Z_{0,up}^{(n)}(i)=Z_{0,down}^{(n)}(i)=0$.

We now define recursively $\tilde{x}_{i}^{(n)}(t)$ for all times $t$. In the
course of the definition, we also define what we mean by "times before the
first jump", "times between up and down jump" and "times between down and up
jump" for a particular process $\tilde{x}_{i}^{(n)}(t)$. We start with

\begin{equation*}
\tilde{x}_{i}^{(n)}(0)=\sqrt{\frac{r(0)}{d^{(n)}(0)}}.  \label{iniialx}
\end{equation*}

In the subcase when $t<\min \{$ $J_{1,down}^{(n)}(i),J_{1,up}^{(n)}(i)\}$ (%
 \textit{times before the first jump}), we have:%
\begin{equation}
\tilde{x}_{i}^{(n)}(t)=\tilde{x}_{i}^{(n)}(0)\exp (-\int_{0}^{t}\frac{b(u)}{2%
}\ud u)+\int_{0}^{t}\exp (-\int_{u}^{t}\frac{b(s)}{2}\ud s)\frac{\sigma (u)}{2}%
\ud W_{i}(u).  \notag
\end{equation}

In the case where $J_{p-1,down}^{(n)}(i)\leq J_{k,up}^{(n)}(i)\leq t<$ $%
J_{p,down}^{(n)}(i)\leq J_{k+1,up}^{(n)}(i)$ (\textit{times between up and
down jump}), we define:%
\begin{equation*}
\tilde{x}_{i}^{(n)}(t)=\tilde{x}_{i}^{(n)}(J_{k,up}^{(n)}(i))(1+\frac{1}{2}%
\frac{\ud d^{(n)}}{\ud t}|_{t=J_{k,up}^{(n)}(i)})\exp
(-\int\limits_{J_{k,up}^{(n)}(i)}^{t}\frac{b(u)}{2}\ud u)+\int%
\limits_{J_{k,up}^{(n)}(i)}^{t}\exp (-\int\limits_{u}^{t}\frac{b(s)}{2}\ud s)%
\frac{\sigma (u)}{2}\ud W_{i}(u).
\end{equation*}

In the case where $J_{k-1,up}^{(n)}(i)\leq J_{p,down}^{(n)}(i)\leq t<$ $%
J_{k,up}^{(n)}(i)\leq J_{p+1,down}^{(n)}(i)$ (\textit{times between down and
up jump}), we define:%
\begin{equation*}
\tilde{x}_{i}^{(n)}(t)=\tilde{x}_{i}^{(n)}(J_{p,down}^{(n)}(i))(1+\frac{1}{2}%
\frac{\ud d^{(n)}}{\ud t}|_{t=J_{p,down}^{(n)}(i)})\exp
(-\int\limits_{J_{p,down}^{(n)}(i)}^{t}\frac{b(u)}{2}\ud u)+\int%
\limits_{J_{p,down}^{(n)}(i)}^{t}\exp (-\int\limits_{u}^{t}\frac{b(s)}{2}\ud s)%
\frac{\sigma (u)}{2}\ud W_{i}(u).
\end{equation*}

We summarize this into the following definitions.

\begin{eqnarray*}
&&\ud \tilde{x}_{i}^{(n)}(t)=-(\frac{b(t)}{2}+\frac{g^{(n)}(t,i)}{2})\tilde{x}%
_{i}^{(n)}(t)\ud t+\frac{\sigma (t)}{2}\ud W_{i}(t),  \label{x_io}\\
&&\mbox{with}\\
&&g^{(n)}(t,i)=\left( \sum\limits_{k=1}^{N_{up}^{(n)}(i,t)}\delta (t-\frac{1}{2%
}J_{k,up}^{(n)}(i))+\sum\limits_{p=1}^{N_{down}^{(n)}(i,t)}\delta
(t-J_{p,down}^{(n)}(i))\right) \frac{\ud d^{(n)}}{\ud t}|_{t}.  \label{g^n}
\end{eqnarray*}

\textbf{Fact 1}

There is no jump of $\tilde{x}_{i}^{(n)}(t)$ in the interval $%
[J_{p,down}^{(n)}(i)\leq t<J_{k(p),up}^{(n)}(i))$. Likewise, there is no
jump of $\tilde{x}_{i}^{(n)}(t)$ in the interval $[J_{k,up}^{(n)}(i)\leq
t<J_{p(k),down}^{(n)}(i))$.

\bigskip

We further define:

\begin{eqnarray*}
W^{(n)}(t,a,\omega ) &=&\sum_{i=1}^{\infty }1\{i-1<d^{(n)}(a)\leq
i\}W_{i}(t,\omega ), \\
x^{(n)}(t,u) &=&\tilde{x}_{i}^{(n)}(t)1\{i-1<u\leq i\},
\end{eqnarray*}
where $1\{i-1<d^{(n)}(a)\leq
i\}$ denotes an indicator function.
\bigskip
We also introduce the following notations:
\begin{eqnarray*}
&&Z_{mid}^{(n,0,M)}(0)\overset{}{=}\int%
\limits_{u=0}^{d^{(n)}(t)}x^{(n)}(0,u)^{2}\ud u, \\
&&\ud\tilde{x}_{i}^{(n,0,M)}(t)=\frac{1}{2}(-b(t)+g^{(n)}(t,i))\tilde{x}%
_{i}^{(n,0,M)}(t)\ud t+\frac{\sigma (t)}{2}\ud W_{i}(t)~\mbox{for}~0\leq t\leq
t_{1}^{M},
\end{eqnarray*}

and for each $1\leq m<M$:

\begin{eqnarray*}
&&\tilde{x}_{i}^{(n,m,M)}(t_{m}^{M})=\sqrt{\frac{%
Z_{mid}^{(n,m,M)}(t_{m}^{M})}{d^{(n)}(t_{m}^{M})}}, \\
&&\ud\tilde{x}_{i}^{(n,m,M)}(t)=\frac{1}{2}(-b(t)+g^{(n)}(t,i))\tilde{x}%
_{i}^{(n,m,M)}(t)\ud t+\frac{\sigma (t)}{2}\ud W_{i}(t)~\mbox{for}~
t_{m}^{M}\leq t\leq t_{m+1}^{M}, \\
&&x_{mid}^{(n,m,M)}(t,u)=\tilde{x}_{i}^{(n,m,M)}(t)1\{i-1<u\leq i\},
\end{eqnarray*}

in which we define $Z_{mid}^{(n,m,M)}(t)$ for all $1\leq m<M$ as:

\begin{equation*}
Z_{mid}^{(n,m,M)}(t)\overset{}{=}\int%
\limits_{u=0}^{d^{(n)}(t)}(x_{mid}^{(n,m,M)}(t,u))^2\ud u~\mbox{for}~t_{m}^{M}\leq
t\leq t_{m+1}^{M}.
\end{equation*}%
\bigskip

We denote by (for the $k$ branches) lower integration bounds:

\begin{equation*}
L_{up}(m,M,k,i,n)=\min (t_{m}^{M},J_{k,up}^{(n)}(i))
\end{equation*}

and upper integration bounds:

\begin{equation*}
U_{up}(m,k,i,n)=\min (J_{next}^{(n)}(J_{k,up}^{(n)}(i)),t_{m+1}^{M}),
\end{equation*}
where $J_{next}^{(n)}$ is defined as in (\ref{Jnext}).

Similarly we can define $L_{down}(m,M,p,i,n)$ and $L_{down}(m,M,p,i,n)$ for the $p$ branches. We define a remainder by

\begin{eqnarray*}
R^{(n,m,M)}(t)=\sum_{i=1}^{M_d^{(n)}}\sum_{k=1}^{N_{up}(i)}\int%
\limits_{L_{up}(m,M,k,i,n)}^{U_{up}(m,M,k,i,n)}\left( (\tilde{x}%
_{i}^{(n,m,M)}(s))^{2}\frac{\ud d^{(n)}}{\ud s}-(\tilde{x}%
_{i}^{(n,m,M)}(J_{k,up}^{(n)}(i)))^{2}\frac{\ud d^{(n)}}{\ud s}%
|_{s=J_{k,up}^{(n)}(i)}\right) \ud s+ \\
\sum_{p=1}^{N_{down}(i)}\int%
\limits_{L_{down}(m,M,p,i,n)}^{U_{down}(m,M,p,i,n)}\left( (\tilde{x}%
_{i}^{(n,m,M)}(s))^{2}\frac{\ud d^{(n)}}{\ud s}-(\tilde{x}%
_{i}^{(n,m,M)}(J_{p,down}^{(n)}(i)))^{2}\frac{\ud d^{(n)}}{\ud s}%
|_{s=J_{p,down}^{(n)}(i)}\right) \ud s.
\end{eqnarray*}

We now define $r_{j,mid}^{(n,m,M)}(t)$ by, for $t_{m}^{M}\leq t<t_{m+1}^{M}$:

\begin{eqnarray*}
&&r_{j,mid}^{(n,m,M)}(t_{m}^{M}) =\frac{Z_{mid}^{(n,m,M)}(t_{m}^{M})}{n}, \\
&&\mbox{and as a consequence,}\\
&&r_{j,mid}^{(n,m,M)}(t)-r_{j,mid}^{(n,m,M)}(t_{m}^{M})\\
&&=-\int_{s=t_{m}^{M}}^{t}(b(s)r_{j,mid}^{(n,m,M)}(s)+\frac{\sigma
^{2}(s)d^{(n)}(s)}{4})\ud s+\int_{s=t_{m}^{M}}^{t}\sigma (s)\sqrt{%
r_{j,mid}^{(n,m,M)}(s)}\ud B_{j,mid}^{(n,m,M)}(s) \\
&&+\frac{R^{(n,m,M)}(t)}{n}.
\end{eqnarray*}

As usual $\{B_{j,mid}^{(n,m,M)}(s)\}_j$ are independent Brownian motions. We can
thus define the process:

\begin{equation}
\label{r1mid}
r_{j,mid}^{(n,M)}(t)=\sum_{m=0}^{M}r_{j,mid}^{(n,m,M)}(t)1{\{t_{m}^{M}\leq
t<t_{m+1}^{M}\}}.
\end{equation}%
\bigskip

\subsubsection{The Scaled Non-central (SNC) Chi-squared Distribution}

For any $\lambda _{1}>0,\lambda _{2}\geq 0$ and $c>0$, and for any $x\in
\mathbb{R}$, we define:%
\begin{equation*}
g_{\lambda _{1},\lambda _{2},c}(x)=\frac{1}{c^{2}}\sum\limits_{i=0}^{\infty }%
\frac{e^{-\lambda _{2}/2}(\lambda _{2}/2)^{i}}{i!}f_{\lambda _{1}+2i}(\frac{x%
}{c^{2}}),
\end{equation*}

where%
\begin{equation*}
f_{\lambda }(x)=\left\{
\begin{array}{c}
\frac{x^{\lambda /2-1}e^{-x/2}}{2^{\lambda /2}\Gamma (\lambda /2)} \\
0%
\end{array}%
\right. \text{if}%
\begin{array}{c}
x\geq 0 \\
\text{otherwise.}%
\end{array}%
\end{equation*}

It is straightforward to verify that $g_{\lambda _{1},\lambda _{2},c}$ is a
probability density. We say that a random variable $X\sim \chi ^{2}(\lambda
_{1},\lambda _{2},c)$ if the density of $X$ is $g_{\lambda _{1},\lambda
_{2},c}$. In words, $X$ is a scaled non-central chi-square with real-valued
degrees of freedom $\lambda _{1}$. When $c=1$ and $\lambda _{1}$ is integer,
we obtain the standard non-central chi-square distribution. When $\lambda
_{1}$ is integer the random variable $X$ is the sum of squares of
independent normal random variables $X_{i}$ with mean $\mu _{i}$ and
variance $c^{2}$:%
\begin{equation*}
X=\sum_{i=1}^{\lambda _{1}}X_{i}^{2}
\end{equation*}

where the parameter
\begin{equation*}
\lambda _{2}=\sum_{i=1}^{\lambda _{1}}\frac{\mu _{i}^{2}}{c^{2}}
\end{equation*}

As the usual chi-square distribution, its generalization to real-valued
degrees of freedom is infinitely divisible.

\bigskip
\begin{theorem}
\label{thm1} For any $X\sim \chi ^{2}(\lambda _{1},\lambda
_{2},c)$ there exists $n$ independent and identically
distributed random variables $X_{1},...,X_{n}$ such that:
\begin{equation*}
X\overset{d}{=}\sum_{k=1}^{n}X_{k}.
\end{equation*}

Moreover:
\begin{equation*}
X_{1}\sim \chi ^{2}(\frac{\lambda _{1}}{n},\frac{\lambda _{2}}{n},c).
\end{equation*}%
\end{theorem}
\bigskip

\subsection{Main Result}
\begin{theorem}
\label{thm2}
Suppose Assumption 1 holds. Then
\begin{description}
\item[(i)] Assume that $d(.)$ is
differentiable over $[0,T]$, bounded away from zero, and has a finite number of optima.
For each $t\in \lbrack 0,T]$ the characteristic
function of $r(t)$ is: for $\omega\in\mathbb R$,

\begin{equation}
\label{master_result}
E[\exp (i\omega r(t))]=\exp \left( i\omega \left( \frac{r(0)e^{-%
\int_0^tb(u)\ud u}}{1-2i\omega \Sigma (0,t)}+\int_{0}^{t}\frac{\exp
(-\int_{s}^{t}b(u)\ud u)\theta (s)}{1-2i\omega \Sigma (s,t)}\ud s\right) \right),
\end{equation}

where

\begin{equation*}
\Sigma (s,t):=\frac{1}{4}\int_{s}^{t}\exp (-\int_{v}^{t}b(u)\ud u)\sigma
^{2}(v)\ud v.
\end{equation*}

\item[(ii)] For any $t\in[0,T]$, the following convergence holds in distribution:
\begin{equation*}
\lim_{n\rightarrow \infty }r_{1,mid}^{(n,n^{4Q})}(t)=r(t),  \label{part2ofth}
\end{equation*}
where the definition of $r_{1,mid}^{(n,n^{4Q})}(t)$ is given in (\ref{r1mid}).
\end{description}
\end{theorem}
\textbf{Remarks:}

\begin{enumerate}
\item The representation (\ref{lemma28}), which together with (\ref{b1}) and (\ref%
{sigma1}), states that:%
\begin{eqnarray*}  \label{beforeIBP}
&& E[\exp (i\omega r(t)]  \notag \\
&&=\frac{\exp \left(\frac{i\omega r(0)e^{-\int_{0}^{t}b(u)\,\mathrm{d}u}}{%
1-2i\omega \int_{0}^{t}e^{-\int_{v}^{t}b(u)\,\mathrm{d}u}\sigma ^{2}(v)\,%
\mathrm{d}v/4}-\frac{1}{2}\int_{0}^{t}d^{\prime }(s)\log \left( 1-2i\omega
\int_{s}^{t}e^{-\int_{v}^{t}b(u)\,\mathrm{d}u}\sigma ^{2}(v)/4\,\mathrm{d}%
v\right) \,\mathrm{d}s\right)}{(1-2i\omega \int_{0}^{t}e^{-\int_{v}^{t}b(u)\,%
\mathrm{d}u}\sigma ^{2}(v)/4\,\mathrm{d}v)^{d(0)/2}}.
\end{eqnarray*}

Therefore Equation (\ref{master_result}) results from the following integration by parts%
\begin{eqnarray*}
&&\frac{1}{2}\int_{0}^{t}d^{\prime }(s)\log \big(1-2i\omega
\int_{s}^{t}e^{-\int_{v}^{t}b(u)\,\mathrm{d}u}\sigma ^{2}(v)/4\,\mathrm{d}t%
\big)\,\mathrm{d}s \\
&=&-\frac{d(0)}{2}\log \big(1-2i\omega \int_{s}^{t}e^{-\int_{v}^{t}b(u)\,%
\mathrm{d}u}\sigma ^{2}(v)/4\,\mathrm{d}t\big)-\frac{1}{2}\Big(\int_{0}^{t}%
\frac{d(s)e^{-\int_{s}^{t}b(u)\,\mathrm{d}u}\sigma ^{2}(s)/4}{(1-2i\omega
\int_{s}^{t}e^{-\int_{v}^{t}b(u)\,\mathrm{d}u}\sigma ^{2}(v)/4\,\mathrm{d}v)}%
\,\mathrm{d}s\Big),
\end{eqnarray*}

\item When $b,\sigma $, and $d$ are constants over $[0,T]$, i.e., $d^{\prime
}(s)=0$, we see from (\ref{beforeIBP}) that $r(t)$ has the SNC chi-squared
distribution with $d(0)$ degrees of freedom. Our characteristic function
thus properly generalizes a well-known result.

\item Casual observation of Lemma \ref{lemma2.8} seems to highlight a much simpler
proof: namely, approximate the time-dependent parameters by step functions,
i.e., construct an Euler approximation of (\ref{CIR}), calculate its characteristic
function like in Lemma \ref{lemma2.8}, and make the time-step go to zero.
Unfortunately, we could not find any proof of strong convergence of the
Euler approximation for the CIR process: in our proof we had to resort to
the extra power coming from our representation (\ref{part2ofth}), and to
weak convergence. Even if there was a classical proof of strong convergence
of the Euler approximation with deterministic volatility, the extension to
stochastic volatility would be difficult, unlike what we obtain in the next
section.

\item We have not found in the literature any distribution with
characteristic function equal to (\ref{master_result}). This seems to be a
new, or at least independently rediscovered distribution. We discuss in the
stochastic volatility section a method to approximate this distribution by
its moments.

\item It can be verified that (\ref{master_result}) solves the equivalent of
the Fokker-Planck equation for characteristic functions. We emphasize again
that the main strength of our method of proof, compared to solving the
Fokker-Planck equation directly, is that it can be extended to stochastic
volatility.
\end{enumerate}

The Fokker-Planck equation for the density $f(r,t)$ of $r(t)$ given in (\ref{CIR}) is:

\begin{eqnarray*}
\frac{\partial f(r,t)}{\partial t} &=&-\frac{\partial }{\partial r}%
[(-b(t)r+\theta(t))f(r,t)]+\frac{1}{2}\frac{\partial ^{2}}{\partial r^{2}}%
[\sigma ^{2}(s)rf(r,t)] \\
&=&b(t)f(r,t)+(\sigma ^{2}(t)+b(t)r-\theta (t))\frac{\partial f(r,t)}{%
\partial r}+\frac{\sigma ^{2}(t)}{2}r\frac{\partial ^{2}f(r,t)}{\partial
r^{2}}.
\end{eqnarray*}

Let the Fourier transform with respect to $x$ be defined as $\hat{f}(x,t)=%
\frac{1}{\sqrt{2\pi }}\int_{} e^{-ixr}f(r,t)\ud r$, then the Fourier transform
with respect to $x$ of the Fokker-Planck equation is given as
\begin{eqnarray}
\frac{\partial \hat{f}(x,t)}{\partial t}&=&-b(t)x\frac{\partial \hat{f}(x,t)%
}{\partial x}-\theta(t) ix\hat{f}(x,t)-\frac{i\sigma ^{2}(t)}{2}x^{2}\frac{%
\partial \hat{f}(x,t)}{\partial x}  \notag \\
&=&-\theta(t) ix\hat{f}(x,t)-(b(t)x+\frac{i\sigma ^{2}(t)}{2}x^{2})\frac{%
\partial \hat{f}(x,t)}{\partial x}.  \label{mainFP}
\end{eqnarray}

Writing $\Phi(x,t)=E[e^{ixr(t)}]$, then:

\begin{equation*}
\Phi(x,t) =\exp\Big(ix\big(\frac{r(0)e^{-\int_0^tb(u)\ud u}}{%
1-2ix\int_0^te^{-\int_v^tb(u)\ud u}\sigma^2(v)/4\ud v}+\int_0^t\frac{%
\theta(s)e^{-\int_s^tb(u)\ud u}}{(1-2ix\int_s^te^{-\int_v^tb(u)\ud
u}\sigma^2(v)/4\ud v)}\ud s\big)\Big).
\end{equation*}
Hence the following relation holds:
\begin{equation*}
\Phi (-x,t) =\sqrt{2\pi }\hat{f}(x,t).
\end{equation*}
We first compute the left hand side (LHS) of (\ref{mainFP}). Observe that

\begin{equation*}
\sqrt{2\pi }\frac{\partial \hat{f}\mathbb{(}x,t\mathbb{)}}{\partial t}%
=\Phi(-x,t)\left(-ix\frac{\partial}{\partial t}\left(\frac{%
r(0)e^{-\int_0^tb(u)\ud u}}{1+2ix\int_0^te^{-\int_v^tb(u)\ud u}\sigma^2(v)/4\ud v}%
+\int_0^t\frac{\theta(s)e^{-\int_s^tb(u)\ud u}}{(1+2ix\int_s^te^{-\int_v^tb(u)%
\ud u}\sigma^2(v)/4\ud v)}\ud s \right)\right).
\end{equation*}

We calculate

\begin{eqnarray*}
&&\frac{\partial}{\partial t}\left( \frac{r(0)e^{-\int_0^tb(u)\ud u}}{%
1+2ix\int_0^te^{-\int_v^tb(u)\ud u}\sigma^2(v)/4\ud v}+\int_0^t\frac{%
\theta(s)e^{-\int_s^tb(u)\ud u}}{(1+2ix\int_s^te^{-\int_v^tb(u)\ud %
u}\sigma^2(v)/4\ud v)}\ud s\right) \\
&&=\frac{r(0)e^{-\int_0^tb(u)\ud u}(-b(t)-ix\sigma^2(t)/2)}{(1+2ix\int_0^te^{-%
\int_v^tb(u)\ud u}\sigma^2(v)/4\ud v)^2}+\theta(t)+\int_{0}^{t}\frac{%
\theta(s)e^{-\int_s^tb(u)\ud u}(-b(t)-ix\sigma^2(t)/2)}{(1+2ix\int_s^te^{-%
\int_v^tb(u)\ud u}\sigma^2(v)/4\ud v)^2}\ud s \\
&&=\theta(t)-\left(b(t)+ix\frac{\sigma^2(t)}{2}\right)\left(\frac{r(0)e^{-%
\int_0^tb(u)\ud u}}{(1+2ix\int_0^te^{-\int_v^tb(u)\ud u}\sigma^2(v)/4\ud v)^2}%
+\int_{0}^{t}\frac{\theta(s)e^{-\int_s^tb(u)\ud u}}{(1+2ix\int_s^te^{-%
\int_v^tb(u)\ud u}\sigma^2(v)/4\ud v)^2}\ud s\right).
\end{eqnarray*}

Secondly we compute the following item for the right hand side (RHS) of (\ref%
{mainFP}):

\begin{equation*}
\sqrt{2\pi }\frac{\partial \hat{f}(x,t) }{\partial x}=\Phi (-x,t)\frac{%
\partial}{\partial x}\left(\frac{-ixr(0)e^{-\int_0^tb(u)\ud u}}{%
1+2ix\int_0^te^{-\int_v^tb(u)\ud u}\sigma^2(v)/4\ud v}+\int_0^t\frac{%
-ix\theta(s)e^{-\int_s^tb(u)\ud u}}{(1+2ix\int_s^te^{-\int_v^tb(u)\ud %
u}\sigma^2(v)/4\ud v)}\ud s \right).
\end{equation*}

On one hand,

\begin{equation*}
\frac{\partial}{\partial x}\left(\frac{-ixr(0)e^{-\int_0^tb(u)\ud u}}{%
1+2ix\int_0^te^{-\int_v^tb(u)\ud u}\sigma^2(v)/4\ud v}\right)= -\frac{%
ir(0)e^{-\int_0^tb(u)\ud u}}{(1+2ix\int_0^te^{-\int_v^tb(u)\ud u}\sigma^2(v)/4\ud
v)^{2}}.
\end{equation*}
On the other hand,
\begin{eqnarray*}
\frac{\partial}{\partial x}\left(\int_0^t\frac{-ix\theta(s)e^{-\int_s^tb(u)\ud
u}}{(1+2ix\int_s^te^{-\int_v^tb(u)\ud u}\sigma^2(v)/4\ud v)}\ud s\right)
=-\int_{0}^{t}\frac{i\theta(s)e^{-\int_s^tb(u)\ud u}}{(1+2ix\int_s^te^{-%
\int_v^tb(u)\ud u}\sigma^2(v)/4\ud v)^{2}}\ud s.
\end{eqnarray*}

Multiplying both sides of Equation (\ref{mainFP}) by $\sqrt{2\pi }/\Phi(-x,t)
$ we get:

\begin{eqnarray*}
&&LHS=-ix\theta(t) \\
&&+ix\left((b(t)+ix\frac{\sigma^2(t)}{2})\left(\frac{%
r(0)e^{-\int_0^tb(u)\ud u}}{(1+2ix\int_0^te^{-\int_v^tb(u)\ud u}\sigma^2(v)/4\ud v)^2%
}+\int_{0}^{t}\frac{\theta(s)e^{-\int_s^tb(u)\ud u}}{(1+2ix\int_s^te^{-%
\int_v^tb(u)\ud u}\sigma^2(v)/4\ud v)^2}\ud s\right)\right).
\end{eqnarray*}

and
\begin{eqnarray*}
&&RHS=-ix\theta(t) \\
&&-(\frac{i\sigma ^{2}(t)}{2}x^{2}+b(t)x)\left( -\frac{%
ir(0)e^{-\int_0^tb(u)\ud u}}{(1+2ix\int_0^te^{-\int_v^tb(u)\ud u}\sigma^2(v)/4\ud
v)^{2}}-\int_{0}^{t}\frac{i\theta(s)e^{-\int_s^tb(u)\ud u}}{(1+2ix\int_s^te^{-%
\int_v^tb(u)\ud u}\sigma^2(v)/4\ud v)^{2}}\ud s\right).
\end{eqnarray*}
It is easy to see LHS and RHS are algebraically equal. Thus the
Fokker-Planck equation (\ref{mainFP}) has been verified.

\section{Option Pricing}

We call $P(t,T)$ the price at time $t$ of a zero-coupon bond with maturity $%
T $, i.e.:%
\begin{equation*}
\label{P}
P(t,T)=E\big[\exp (-\int_{t}^{T}r(s)\ud s)|\mathcal{F}_{t}\big].
\end{equation*}%
\bigskip

A standard result (see, e.g. Shreve (2004)) is that:%
\begin{equation*}
\label{CA}
P(t,T)=\exp \big(-r(t)C(t,T)+A(t,T)\big),
\end{equation*}

where $C(t,T)$ satisfies the Riccati equation :%
\begin{eqnarray}
C_{t}(t,T)-b(t)C(t,T)-\frac{1}{2}\sigma ^{2}(t)C(t,T)^{2}+1 &=&0
\label{equ_forb} \\
C(T,T) &=&0  \notag
\end{eqnarray}

and:%
\begin{equation}
A(t,T)=-\int_{t}^{T}\theta (u)C(u,T)\ud u.  \label{equforA}
\end{equation}

There is a well-known analytical solution in case the parameters $b,\sigma $
and $\theta $ are constant.

\subsection{Call Option Price}

By moving to the forward measure, we can separate discounting and evaluation
of the expected value of the payoff. We can thus apply Theorem \ref{thm2}.

\bigskip

Let now $0\leq t\leq T$. The price $C(0)$ at initial time $0$ of a European option
(expiring at time $t$ and with exercise price $K$)\ on the $T$ maturity bond
is given by:%
\begin{equation*}
C(0)=E\Big[\exp (-\int\limits_{0}^{t}r(u)\ud u)\max (P(t,T)-K,0)\Big].
\end{equation*}

The $t$-forward probability measure is the probability measure $\mathbb{P}%
^{t}$ under which $P(.,t)$ is a martingale. In this measure the process $
\{B^{t}(s)\}_s$ is Brownian motion where:%
\begin{equation*}
\ud B^{t}(s)=\ud B(t)-\sigma _{P}(s,t)\ud t,
\end{equation*}%
\bigskip

where the bond volatility $\sigma _{P}(s,t)$ is given by:%
\begin{equation*}
\sigma _{P}(s,t)=\left\{
\begin{array}{c}
-\sigma (s)\sqrt{r(s)}C(s,t) \\
0%
\end{array}%
\right. \text{if}%
\begin{array}{c}
0\leq s\leq t \\
\text{else.}%
\end{array}%
\end{equation*}

Thus:%
\begin{equation*}
\ud r(s)=((-b(s)+\sigma ^{2}(s)C(s,t))r(s)+\theta (s))\ud s+\sigma (s)\sqrt{r(s)}%
\ud B^{t}(s).
\end{equation*}

We define:%
\begin{eqnarray*}
b^{t}(s) &=&-b(s)+\sigma ^{2}(s)C(s,t) \\
\Sigma ^{t}(s,t) &=&\frac{1}{4}\int_{0}^{t}\exp
(-\int_{s}^{t}b^{t}(u)\ud u)\sigma ^{2}(s)\ud s
\end{eqnarray*}

as well as the Laplace transform of the density of $r(t)$ in the $t$-forward
measure:

\begin{equation*}
\hat{F}_{r(t)}(p)=\exp \left( -p\left( \frac{r(0)e^{-\int_0^tb^t(u)\ud u}}{%
1+2p\Sigma ^{t}(0,t)}+\int_{0}^{t}\frac{\exp (-\int_{s}^{t}b^{t}(u)\ud u)\theta
(s)}{1+2p\Sigma ^{t}(s,t)}\,\mathrm{d}s\right) \right).
\end{equation*}
\begin{theorem}
\label{thm3}\textit{\ } \textit{Let }$C(s,T)$\textit{\ and }$A(s,T)$%
\textit{\ solve (\ref{equ_forb})\ and} (\ref{equforA}). \textit{The price of
a call option\ is given by Laplace inversion:}

\begin{equation}
\label{C0}
C(0)=\frac{P(0,t)}{2\pi i}\lim\limits_{b\rightarrow \infty
}\int\limits_{p=a-ib}^{a+ib}e^{pr_{0}}\hat{C}_{t}(p;t)\hat{F}_{r(t)}(p)\ud p,
\end{equation}

where,

\begin{equation}
\label{Cp}
\hat{C}_{t}(p)=\frac{pK^{p-C(t,T)}e^{A(t,T)}-(p+C(t,T))K^{p+1}}{C(t,T)p-p^2
}.
\end{equation}
\end{theorem}
\textbf{Proof }Using the change of measure developed by Jamshidian (1989)
and Geman (1989):%
\begin{equation*}
C(0)=P(0,t)E^{\mathbb{P}^{t}}\Big[\max (\exp (-r(t)C(t,T)+A(t,T))-K,0\Big].
\end{equation*}

The Laplace transform of $\max (P(t,T)-K,0)$ is (see, e.g. Lewis\ (2000)):%
\begin{equation*}
\int\limits_{0}^{\infty }e^{-pt}\max \Big(e^{-C(t,T)r+A(t,T)}-K,0\Big)\ud r=\hat{C}%
_{t}(p),
\end{equation*}
which is equal to the right hand side of Equation (\ref{Cp}). 
It is well-known (see, e.g., Lewis (2000) again) that the price of an option
is given by the discounted inverse Laplace transform of the product of the
transform of the payoff and the transform of the distribution. Therefore we obtain (\ref{C0}). $\square$

\section{Stochastic Volatility}

We now assume that both $\sigma $ and $\theta $ are continuous functions of
infinite variation, but that the dimension $d(t)$ satisfies assumption 1.
This assumption is not unrealistic, in the sense that only the moments of $%
\theta $ intervene in the calculations of the moments of $r(t)$. Given a
time-series of $r$, it would actually be difficult to invalidate the
hypothesis that $\theta $ is of infinite variation. Thus, given the limited
practical use of generalizing our results to the case where $d$ can be of
infinite variation, we leave this task as an open conjecture for future
research. There is an interesting parallel between the deterministic and the
stochastic volatility case. In the deterministic case it was deemed natural
since Magshoodi (1996) to take the simplifying assumption that $d$ be
integer-valued. In the stochastic case, we deem natural to consider the case
where $d$ if of finite variation.

We now consider two cases. The easiest case to consider is when the
volatility is independent on the interest rate. We then build a tractable
model where volatility is correlated with the interest rate.

\subsection{Independent Volatility}

We consider the model

\begin{eqnarray}
\ud r(t) &=&(-b(t)r(t)+\theta (t))\ud t+\sqrt{v(t)r(t)}\ud B(t),\nonumber\\
r(0)&=&r_{0}>0;  \\
\ud v(t) &=&\mu (v(t),t)\ud t+\xi (v(t),t)\ud B^{v}(t),\nonumber\\
v(0)&=&v_{0}>0,  \label{ind_vol}
\end{eqnarray}

where the coefficients $\mu $ and $\xi $ are such that $0<v(t)<K$ almost
surely. We assume
\begin{equation*}
\ud B(t)\ud B^{v}(t)=0.
\end{equation*}

Under Assumption 1, all the arguments in the proof of Theorem \ref{thm2} hold, and, by conditioning on
the path of volatility, we obtain the characteristic function of the rate:

\begin{equation}
E[\exp (i\omega r(t))]=E\left[\exp \left( i\omega \left( \frac{%
r(0)e^{-\int_0^tb(u)\ud u}}{1-2i\omega \Sigma (0,t)}+\int_{0}^{t}\frac{\exp
(-\int_{s}^{t}b(u)\ud u)\theta (s)}{1-2i\omega \Sigma (s,t)}\,\mathrm{d}%
s\right) \right) \right]
\end{equation}

When volatility is stochastic, it is often possible to write the price of a
bond at time $t$ as a function $f$ of both the rate and volatility:%
\begin{equation*}
P(t,T,\omega )=f(r(t,\omega ),v(t,\omega ),t;T)
\end{equation*}

In this case, to calculate the price $C(0)$ of a bond option, one is more
interested in the conditional Laplace transform function of the density of
the rate, given volatility $v(t)$ at time $t$, which we call $\hat{F}%
_{r(t)|v(t)}$. One can then develop the latter in a Taylor series:
\begin{eqnarray*}
\hat{F}_{r(t)|v(t)}(p) &=&E\left[\exp \left( -p\left( \frac{r(0)e^{-%
\int_0^tb^t(u)\ud u}}{1+2p\Sigma ^{t}(0,t)}+\int_{0}^{t}\frac{\exp
(-\int_{s}^{t}b^{t}(u)\ud u)\theta (s)}{1+2p\Sigma ^{t}(s,t)}\,\mathrm{d}%
s\right) \right) \Big|v(t)\right] \\
&=&\sum_{k=0}^{\infty }(\frac{-p}{k!})^{k}\frac{\ud^{k}\hat{F}_{r(t)|v(t)}}{%
\ud p^{k}}\Big|_{p=0}.
\end{eqnarray*}

For instance the first term is:%
\begin{equation}
\frac{\ud\hat{F}_{r(t)|v(t)}(p)}{\ud p}\Big|_{p=0}=-E\left[(r_{0}e^{-\int_0^tb^t(u)\ud u}+%
\int_{0}^{t}\exp (-\int_{s}^{t}b^{t}(u)\ud u)\theta (s))\Big|v(t)\right].
\label{first_moment}
\end{equation}

It is well-known that the derivative of the Laplace transform of the density
is equal to the negative of the first moment. By It\^o's lemma, we can verify
that the (conditional)\ first moment of the rate is equal to the negative of
the right hand-side of (\ref{first_moment}), which provides a pleasant
confirmation of our result. Since $\theta (s)=d(s)\sigma ^{2}(s)/4$, it is
enough to compute of $E[\sigma ^{2}(s)|v(t)]$ in order to compute explicitly
(\ref{first_moment}). We can then easily reconstruct the conditional
distribution of $r(t)$ given $v(t)$ by inverse Laplace transformation,
provided convergence conditions are met.

\subsection{A Tractable Model with Correlated Volatility}

We consider the following model of volatility, which is a special case of (%
\ref{ind_vol}).%
\begin{eqnarray}
&&\ud w(t ) =(-b(t)w(t )+\theta _{w}(t))\ud t+\sqrt{v(t
)w(t )}\ud B(t ),\nonumber\\
&&r(0)=r_{0}>0; \nonumber\\
&&\ud v(t ) =(-b(t)v(t )+\theta _{v}(t))\ud t+\xi v(t
)\ud B^{v}(t ),\nonumber\\
&&v(0)=v_{0}>0; \nonumber\\
&&r(t ) =w(t)+v(t ).  \label{r_dependent}
\end{eqnarray}

We let%
\begin{eqnarray*}
d_{w}(t) &=&\frac{4\theta _{w}(t)}{v^{2}(t)}, \\
d_{v}(t) &=&\frac{\theta _{w}(t)}{4v^{2}(t)}
\end{eqnarray*}

and suppose that both of them satisfy assumption (1). Even if $B$ and $B^{v}$
are uncorrelated, $r(t)$ is correlated with $v(t)$. It is then clear that

\begin{eqnarray*}
E[\exp (i\omega r(t))|v(t)] &=&E\left[\exp \left( i\omega \left( \frac{%
w(0)e^{-\int_0^tb(u)\ud u}}{1-2i\omega \Sigma (0,t)}+\int_{0}^{t}\frac{%
\exp (-\int_{s}^{t}b(u)\ud u)\theta (s)}{1-2i\omega \Sigma (s,t)}\,\mathrm{d}%
s\right) \right) \Big|v(t)\right] \\
&&\times\exp (i\omega v(t)).
\end{eqnarray*}

The advantage of this model is that the conditional characteristic function
of the rate is explicit. Is this model a stochastic CIR model? We would like
to leave this exciting (and, as far as we know, unsolved) question as an
open problem to the community.

\bigskip

\textbf{Open problem:}

Is there a Brownian motion $B_{r}(t)$ such that the distribution of $r(t)$
defined in (\ref{r_dependent}) agrees with the distribution at time $t$ of
the solution of the following SDE:%
\begin{equation}
\ud x(t)=-b(t)x(t)\ud t+(\theta _{w}(t)+\theta _{v}(t))\ud t+\sqrt{v(t)x(t)}\ud B_{r}(t)%
\text{ }?  \label{dx}
\end{equation}

\textbf{Observation}

One avenue of proof is to reuse our framework. Indeed, let:%
\begin{eqnarray*}
\ud\tilde{x}_{i,w}^{(n)}(t) &=&-(\frac{b(t)}{2}+\frac{g_{w}^{(n)}(t,i)}{2})%
\tilde{x}_{i}^{(n)}(t)\ud t+\frac{\sqrt{v(t)}}{2}\ud W_{i,w}(t), \\
\ud\tilde{x}_{i,v}^{(n)}(t) &=&-(\frac{b(t)}{2}+\frac{g_{v}^{(n)}(t,i)}{2})%
\tilde{x}_{i}^{(n)}(t)\ud t+\frac{\sqrt{v(t)}}{2}\ud W_{i,v}(t),
\end{eqnarray*}%
where $W_{i,w}$ and $W_{i,v}$ are independent Brownian motions, and $%
g_{w}^{(n)}(t,i)$ and $g_{v}^{(n)}(t,i)$ are defined in a manner similar
with (\ref{g^n}). Then, (conditionally on $v$), independent copies $%
w_{j,low}^{(n,M)}$ and $v_{j,low}^{(n,M)}$ exist that satisfy the
representation:%
\begin{eqnarray*}
\sum_{j=1}^{n}w_{j,low}^{(n,M)}(t_{m}^{M})
&=&\sum_{i=1}^{[nd_{w}(t)]}x_{i,w}^{2}(t_{m}^{M}) \\
\sum_{j=1}^{n}v_{j,low}^{(n,M)}(t_{m}^{M})
&=&\sum_{i=1}^{[nd_{v}(t)]}x_{i,v}^{2}(t_{m}^{M})
\end{eqnarray*}

It then remains to prove that, for some appropriately defined processes
$w^{(n,M)}_{1,mid}$ and $v^{(n,M)}_{1,mid}$:
\begin{equation*}
r^{(n,M)}\equiv w^{(n,M)}_{1,mid}+v^{(n,M)}_{1,mid}
\end{equation*}

converges somehow to a process $r^{(\infty ,\infty )}$ that satisfies (\ref%
{dx}), and for which $r^{(\infty ,\infty )}(t)$ agrees with $r(t)$ in
distribution.
\bigskip

\textbf{Acknowledgements}

We would like to thank the participants at the 2014 Claremont Symposium on
Interest rates, as well as the participants of the 2014 Bachelier
conference. We would also thank Professor John Angus for the valuable communication with him. We received very helpful advice from two anonymous referees. All
errors are ours.

\section{Appendix}

\subsection{Proof of Theorem \ref{thm1}}

The characteristic function of $Y\sim \chi ^{2}(\lambda _{1})$ when $\lambda
_{1}>0$ is real-valued is:
\begin{equation}
\Phi _{Y}(\omega )=(1-2i\omega )^{-\lambda _{1}/2}.  \label{ch1}
\end{equation}%
Let $X\sim \chi ^{2}(\lambda _{1},\lambda _{2},1)$. It follows from (\ref%
{ch1}) that:
\begin{eqnarray}
\Phi _{X}(\omega ) &=&\int_{\mathbb{R}}e^{i\omega x}g_{\lambda _{1},\lambda
_{2}}(x)\,\mathrm{d}x=\sum_{k=0}^{+\infty }\frac{e^{-\lambda _{2}/2}(\lambda
_{2}/2)^{k}}{k!}\int_{0}^{+\infty }e^{i\omega x}f_{\lambda _{1}+2k}(x)\,%
\mathrm{d}x  \notag \\
&=&\sum_{k=0}^{+\infty }\frac{e^{-\lambda _{2}/2}(\lambda _{2}/2)^{k}}{k!}%
(1-2i\omega )^{-\lambda _{1}/2-k}=(1-2i\omega )^{-\lambda
_{1}/2}\sum_{k=0}^{+\infty }\frac{e^{-\lambda _{2}/2}((\lambda
_{2}/2)(1-2i\omega )^{-1})^{k}}{k!}  \notag \\
&=&\frac{\exp (\frac{i\lambda _{2}\omega }{1-2i\omega })}{(1-2i\omega
)^{\lambda _{1}/2}}.  \label{charofX2}
\end{eqnarray}%
\bigskip On the other hand, since $X_{1},\ldots ,X_{n}$ are i.i.d,
\begin{equation*}
\Phi _{\sum_{k=1}^{n}X_{k}}(t)=(\Phi _{X_{1}}(t))^{n}.
\end{equation*}%
By the facts that $X=\sum_{k=1}^{n}X_{k}$, we obtain
\begin{equation*}
\Phi _{X_{1}}(\omega )=(\Phi _{X}(\omega ))^{1/n}=\frac{\exp (\frac{%
i(\lambda _{2}/n)\omega }{1-2i\omega })\exp (\frac{2ki\pi }{n})}{%
(1-2it)^{(\lambda _{1}/n)/2}},~\mbox{for $k=0,\ldots,n-1$}.
\end{equation*}%
The fact that $\Phi _{X_{1}}(0)=1$ guarantees
\begin{equation*}
\Phi _{X_{1}}(\omega )=\frac{\exp (\frac{i(\lambda _{2}/n)\omega }{%
1-2i\omega })}{(1-2i\omega )^{(\lambda _{1}/n)/2}}.
\end{equation*}%
Again, it results from (\ref{charofX2}) and the above equation that:
\begin{equation*}
X_{1}\sim \chi ^{2}\big(\frac{\lambda _{1}}{n},\frac{\lambda _{2}}{n},1\big).
\end{equation*}%
The case where $c\neq 1$ follows trivially. $\square$

\subsection{Proof of Theorem \ref{thm2}}

The proofs of the lemmas necessary to this proof are in a second section of Appendix.

\bigskip

For convenience, we define a function $J_{next}^{(n)}(t)$ which associates to a point
in time $t$ the next jump after $t$: 

\begin{equation}
\label{Jnext}
J_{next}^{(n)}(t)=\min \{J\in \mathcal{J}^{(n)}|J>t).
\end{equation}

\bigskip

The stochastic Fubini theorem (Ikeda and Watanabe 1981, Williams 1979 p. 44,
Heath Jarrow Morton 1992, Lemma 0.1) states the following.
\begin{lemma}
\label{lemma0.1}[Heath-Jarrow-Morton (1992)]

Let $\{\Phi (t,a,\omega ):(t,a)\in \lbrack 0,\tau ]\times \lbrack 0,\tau ])$
be a family of real random variables such that:

(i)\ $((t,\omega ),a)\in \{([0,\tau ]\times \Omega )\times \lbrack 0,\tau
]\}\rightarrow \Phi (t,a,\omega )$ is $L\times B[0,\tau ]$ measurable where $%
L$ is the predictable $\sigma $-field;

(ii) $\int\limits_{0}^{t}\Phi ^{2}(s,a,\omega )\ud s<\infty $ \ a.e. for all $%
t\in \lbrack 0,\tau ]$ and $\omega\in\Omega$;

(iii) $\int\limits_{0}^{t}\left( \int_{0}^{\tau }\Phi (s,a,\omega )\ud a\right)
^{2}\ud s<\infty $ \ a.e. for all $t\in \lbrack 0,\tau ]$ and $\omega\in\Omega$.

If $t\rightarrow \int_{0}^{\tau }\int_{0}^{t}\Phi (s,a,\omega )\ud W_{s}(\omega)\ud a$ is
continuous a.e. for all $%
t\in \lbrack 0,\tau ]$ and $\omega\in\Omega$, then:%
\begin{equation*}
\int_{0}^{t}\int_{0}^{\tau }\Phi (s,a,\omega )\ud a\ud W_{s}(\omega)=\int_{0}^{\tau
}\int_{0}^{t}\Phi (s,a,\omega )\ud W_{s}(\omega)\ud a,
\end{equation*}
a.e. for all $t\in \lbrack 0,\tau ]$ and $\omega\in\Omega$.
\end{lemma}
\begin{definition}
 We define:%
\begin{eqnarray*}
\mu ^{(n)}(t,u ) &=&\sum_{i=1}^{\infty }(-(b(s)+g(t,i))\tilde{x}_{i}^{(n)}(t
)^{2}-\frac{\sigma ^{2}(t)}{4})1\{i-1<u\leq i\}; \\
\Phi ^{(n)}(t,u ) &=&\sum_{i=1}^{\infty }\sigma (t)x^{(n)}(t,u); \\
y^{(n)}(t,u) &=&x^{(n)}(t,u)^{2}; \\
Z_{mid}^{(n)}(t) &=&\int\limits_{u=0}^{d^{(n)}(t)}y^{(n)}(t,u)\ud u.
\end{eqnarray*}
\end{definition}
The SDE for $y^{(n)}(t,u)$ spells then: 
\begin{equation*}
\ud y^{(n)}(t,u )=\mu ^{(n)}(t,u )\ud t+\Phi ^{(n)}(t,u
)\partial_{t} W^{(n)}(t,u).
\end{equation*}

Using the convention that the differential operator $d_{s}$ applies to the
first parameter of $W^{(n)}(s,a)$, this definition results in:%
\begin{equation*}
\int_{0}^{\tau }\int_{0}^{t}\Phi ^{(n)}(s,a
)\partial_{s}W^{(n)}(s,a)\ud a=\sum_{i=1}^{\infty }\int_{0}^{\tau
}1\{i-1<d^{(n)}(a)\leq i\}\int_{0}^{t}\sigma (s)x_{i}^{(n)}(s)\ud W_{i}(s)\ud a.
\end{equation*}

We have then the obvious lemma.

\bigskip
\begin{lemma}
\label{lemma2.1} For a.e. $t\in[0,\tau]$, the following equation holds almost surely:
\begin{equation*}
\int_{s=0}^{t}\int_{a=0}^{\tau }\Phi ^{(n)}(s,a
)\ud a\partial_{s}W^{(n)}(s,a)=\int_{a=0}^{\tau }\int_{s=0}^{t}\Phi^{(n)} (s,a
)\partial_{s}W(s,a)\ud a.
\end{equation*}
\end{lemma}
\bigskip

\begin{definition} We define the remainder:

\begin{eqnarray*}
&&R^{(n)}=\sum_{i=1}^{M_d^{(n)}}\sum_{k=1}^{N_{up}(i)}\int%
\limits_{s=J_{k,up}^{(n)}(i)}^{J_{next}^{(n)}(J_{k,up}^{(n)}(i))}\left( (%
\tilde{x}_{i}^{(n)}(s))^{2}\frac{\ud d^{(n)}}{\ud s}|_{s}-(\tilde{x}%
_{i}^{(n)}(J_{k,up}^{(n)}(i)))^{2}\frac{\ud d^{(n)}}{\ud s}|_{J_{k,up}^{(n)}(i)}%
\right) \ud s+ \\
&&\sum_{p=1}^{N_{down}(i)}\int%
\limits_{s=J_{p,down}^{(n)}(i)}^{J_{next}^{(n)}(J_{k,up}^{(n)}(i))}\left( (%
\tilde{x}_{i}^{(n)}(s))^{2}\frac{\ud d^{(n)}}{\ud s}|_{s}-(\tilde{x}%
_{i}^{(n)}(J_{p,down}^{(n)}(i)))^{2}\frac{\ud d^{(n)}}{\ud s}%
|_{J_{p,down}^{(n)}(i)}\right) \frac{\ud d^{(n)}}{\ud s}|_{s}\ud s.
\end{eqnarray*}
\end{definition}
\begin{lemma}
\label{lemma2.2} \textit{There exists a collection of Brownian
motions} $B_{Z,mid}^{(n)}$ \textit{such that:}%
\begin{equation*}
Z_{mid}^{(n)}(t)-Z_{mid}^{(n)}(0)=\int_{s=0}^{t}-b(s)Z_{mid}^{(n)}(s)-\frac{%
\sigma ^{2}(s)d^{(n)}(s)}{4}\ud s+\int_{s=0}^{t}\sigma (s)\sqrt{Z_{mid}^{(n)}(s)%
}dB_{Z,mid}^{(n)}(s)+R^{(n)}(t).  \label{Znmideq}
\end{equation*}
\end{lemma}
Fix a time $t_{m}^{M}$. We cannot estimate the remainder $R^{(n)}(t)$
properly since each term $\tilde{x}_{i}^{(n)}$ may have jumped over $%
[0,t_{m}^{M}]$, which makes them of different orders of magnitude. To this
effect, we refine our representation by resetting it at each time $t_{m}^{M}$
. This allows us to better control the remainder.
\begin{lemma}
\label{lemma2.3} \textit{There exists a collection of Brownian
motions} $B_{Z,mid}^{(n,m,M)}$ \textit{such that, for} $t_{m}^{M}\leq t\leq
t_{m+1}^{M}$:%
\begin{eqnarray*}
&&Z_{mid}^{(n,m,M)}(t)-Z_{mid}^{(n,m,M)}(t_{m}^{M})
=\int_{s=t_{m}^{M}}^{t}-b(s)Z_{mid}^{(n,m,M)}(s)-\frac{\sigma
^{2}(s)d^{(n)}(s)}{4}\ud s \\
&&+\int_{s=t_{m}^{M}}^{t}\sigma (s)\sqrt{Z_{mid}^{(n,m,M)}(s)}%
\ud B_{Z,mid}^{(n,m,M)}(s)+R^{(n,m,M)}(t).
\end{eqnarray*}
\end{lemma}
We notice that the proof of Lemma \ref{lemma2.2} can be carried over exactly to Lemma \ref{lemma2.3}, so we omit it.

\bigskip
\begin{lemma}
\label{lemma2.4} Let $M(n)=n^{4Q}$ .The remainder $%
R^{(n,m,M(n))}(t)/n$ converges to zero in probability when $n\rightarrow
\infty $.
\end{lemma}
\bigskip

We define

\begin{equation*}
\tilde{Z}_{mid}^{(n,M)}(t)=\sum_{m=0}^{M}Z_{mid}^{(n,m,M)}(t)1\{t_{m}^{M}%
\leq t<t_{m+1}^{M}\}
\end{equation*}

By the fact that $Z_{mid}^{(n,m-1,M)}(t_{m}^{M})=Z_{mid}^{(n,m,M)}(t_{m}^{M})$, the trajectory $\tilde{Z}%
_{mid}^{(n,M)}(.,\omega )$ is continuous.\bigskip
\begin{lemma}
\label{lemma2.5} The trajectory\textbf{\ }$%
r_{j,mid}^{(n,m,M)}(.,\omega )$ is continuous almost surely.  Also, the
following relations hold:
\begin{eqnarray*}
r_{j,mid}^{(n,m-1,M)}(t_{m}^{M}) &=&r_{j,mid}^{(n,m,M)}(t_{m}^{M}) \\
\tilde{Z}_{mid}^{(n,M)}(t_{m+1}^{M})
&=&\sum_{j=1}^{n}r_{j,mid}^{(n,m,M)}(t_{m+1}^{M}).  \label{infinite_div}
\end{eqnarray*}
\end{lemma}
\begin{lemma}
\label{lemma2.6} Take $M(n)=n^{4Q}$. The process $
r_{j,mid}^{(n,M(n))}$ converges to $r$ in distribution when $n$ goes to $%
\infty $.
\end{lemma}
\bigskip

\begin{definition}

We define:

\begin{eqnarray*}
d_{low}^{{}}(s,u) &=&[\min_{s\leq t<u}d^{(n)}(t)]; \\
Z_{\chi ^{2}}^{(n,m,M)}(t) &=&\sum_{i=1}^{d_{low}(t_{m}^{M},t_{m+1}^{M})}(%
\tilde{x}_{i}^{(n,m,M)}(t))^{2}\text{ \ \ }t_{m}^{M}\leq t<t_{m+1}^{M};
\label{Z-chideux} \\
\tilde{Z}_{\chi ^{2}}^{(n,M)}(t) &=&\sum_{m=1}^{M_d^{(n)}} Z_{\chi
^{2}}^{(n,m,M)}(t)1\{t_{m}^{M}\leq t<t_{m+1}^{M}\}.
\end{eqnarray*}
\end{definition}
We also say that a $\mathcal{F}_{t}$ - measurable random variable $X$
satisfies

\begin{equation*}
X=o(\Delta ;t)
\end{equation*}

if%
\begin{equation*}
E[X(t)]=o(\Delta ).
\end{equation*}

Each term in the sum in the RHS\ of (\ref{Z-chideux}) has (conditionally on $%
\mathcal{F}_{t_{m}^{M}}$) the SNC chi-squared distribution with one degree
of freedom, and the same scale factor $c^{(m,n,M)}$. Thus $Z_{\chi
^{2}}^{(n)}(t_{m}^{M})$ has (conditionally on $\mathcal{F}_{t_{m}^{M}}$) the
SNC chi-squared distribution with $d_{low}(t_{m}^{M},t_{m+1}^{M})$ degrees
of freedom, and scale factor $c^{(m,n,M)}$. By the infinite divisibility of
the SNC chi-squared distribution (see Theorem \ref{thm1}), we can define $j$
conditionally independent copies $r_{j,low}^{(n,m,M)}$ such that
\begin{equation}
\sum_{j=1}^{n}r_{j,low}^{(n,m,M)}(t)=Z_{\chi ^{2}}^{(n,m,M)}(t)
\label{rjlowchisquare}
\end{equation}%
\bigskip
\begin{lemma}
\label{lemma2.7}
There exists continuous functions $\tilde{b}$ and $\tilde{\sigma}$ such that:%
\begin{eqnarray*}
&&E[\exp (i\omega \tilde{Z}_{mid}^{(n,M)}(t)/n]\\
&&=\frac{\exp (\frac{i\omega \tilde{Z}_{mid}^{(n)}(0)e^{-\int_{0}^{t}\tilde{b}%
(u)\,\mathrm{d}u}}{1-2i\omega \int_{0}^{t}e^{-\int_{v}^{t}\tilde{b}(u)\,%
\mathrm{d}u}\tilde{\sigma}^{2}(v)/4\,\mathrm{d}v}-\frac{1}{2}%
\int_{0}^{t}d^{\prime }(s)\log \left( 1-2i\omega \int_{s}^{t}e^{-\int_{v}^{t}%
\tilde{b}(u)\,\mathrm{d}u}\tilde{\sigma}^{2}(v)/4\,\mathrm{d}v\right) \,%
\mathrm{d}s)}{(1-2i\omega \int_{0}^{t}e^{-\int_{v}^{t}\tilde{b}(u)\,\mathrm{d%
}u}\tilde{\sigma}^{2}(v)/4\,\mathrm{d}v)^{d(0)/2}}+o(\frac{n^{4Q}}{M}).  \notag
\end{eqnarray*}
\end{lemma}
\begin{lemma}
\label{lemma2.8} Take $M(n)=n^{4Q}$. There exists continuous functions $%
\tilde{b}$ and $\tilde{\sigma}$ such that
\begin{eqnarray}
\label{lemma28} 
&&\lim_{n\rightarrow \infty }E[\exp (i\omega r_{mid}^{(n,M(n))}(t)]\nonumber\\
&&=\frac{\exp \left(\frac{i\omega r(0)e^{-\int_{0}^{t}\tilde{b}(u)\,\mathrm{d}u}}{%
1-2i\omega \int_{0}^{t}e^{-\int_{v}^{t}\tilde{b}(u)\,\mathrm{d}u}\tilde{%
\sigma}^{2}(v)/4\,\mathrm{d}v}-\frac{1}{2}\int_{0}^{t}d^{\prime }(s)\log
\left( 1-2i\omega \int_{s}^{t}e^{-\int_{v}^{t}\tilde{b}(u)\,\mathrm{d}u}%
\tilde{\sigma}^{2}(v)/4\,\mathrm{d}v\right) \,\mathrm{d}s\right)}{(1-2i\omega \int_{0}^{t}e^{-\int_{v}^{t}\tilde{b}(u)\,\mathrm{d%
}u}\tilde{\sigma}^{2}(v)/4\,\mathrm{d}v)^{d(0)/2}}. 
\end{eqnarray}
\end{lemma}
Taking the sequence $M(n)=n^{4Q}$, Lemma \ref{lemma2.6} and Lemma \ref{lemma2.8} show that $r(t)$ has
characteristic function given by (\ref{lemma28}). Formula (\ref{master_result}) obtains by integration by parts. Since both forms of the
characteristic functions (before and after integration by parts) are
interesting, we show this calculation in the main text. We calculate the
first two moments of $r(t)$ by ways of the characteristic function (\ref{lemma28}) and match them with the analytical results obtained from
applying It\^o's lemma to (\ref{masterstoch2}). Using (\ref{master_result}),
it is easy to see that the only valid choice (for all $t$)\ of a function $%
\tilde{b}$ and $\tilde{\sigma}$ is:%
\begin{eqnarray}
\tilde{b} &=&b  \label{b1} \\
\tilde{\sigma} &=&\sigma  \label{sigma1}
\end{eqnarray}

Note that it is also possible to see directly from our definitions that (\ref%
{b1}), (\ref{sigma1}) hold, but it is more cumbersome.

\subsection{Proofs of All Lemmas}

\noindent \textbf{Proof of Lemma \ref{lemma2.2}}

The function $t\rightarrow \int_{0}^{\tau }\int_{0}^{t}\Phi
^{(n)}(s,a )\ud W_{s}(a)\ud a$ is continuous a.e. and a.s. even at a discontinuity
point\footnote{%
The function $t\rightarrow \int_{0}^{t}\Phi ^{(n)}(s,a )\ud W_{s}$ is
not continuous at a discontinuity $t=d^{-1}(i/n)$, but the discontinuity is
"integrated" when this function is integrated w.r.t. $a$}. We note however
that the function:%
\begin{equation*}
t\rightarrow \int\limits_{u=0}^{d^{(n)}(t)}\int\limits_{s=0}^{t}\mu
^{(n)}(s,a,\omega )\ud s\ud a
\end{equation*}

is discontinuous at $t$. The latter case is however covered by the regular
Fubini theorem (see e.g. Hunter and Nachtergaele 2001, p.350). Applying
both Fubini theorems we have:%
\begin{equation}
\label{ZABC}
Z_{mid}^{(n)}(t)-Z_{mid}^{(n)}(0)=(A)+(B)+(C),
\end{equation}

where:%
\begin{eqnarray}
(A)&=&\int\limits_{s=0}^{t}\int\limits_{u=0}^{d^{(n)}(s)}\mu ^{(n)}(s,u)\ud u%
\ud s;  \notag \\
(B)&=&\int\limits_{s=0}^{t}\int\limits_{u=0}^{d^{(n)}(s)}\Phi
^{(n)}(s,u)\ud u\partial_{s}W(s,u);  \notag \\
(C)&=&\int\limits_{s=0}^{t}\int\limits_{u=d^{(n)}(s)_{+}}^{d^{(n)}(t)}\mu
^{(n)}(s,u)\ud u\text{ }\ud s+\int\limits_{s=0}^{t}\int%
\limits_{u=d^{(n)}(s)_{+}}^{d^{(n)}(t)}\Phi ^{(n)}(s,u)\ud u\partial_{s}W(s,u)  \notag
\\
&=&\int\limits_{s=0}^{t}\int\limits_{a=s_{+}}^{t}\mu ^{(n)}(s,d^{(n)}(a))%
\frac{\ud d^{(n)}}{\ud a}|_{a}\ud a\ud s+\int\limits_{s=0}^{t}\int\limits_{a=s_{+}}^{t}%
\Phi ^{(n)}(s,d^{(n)}(a))\frac{\ud d^{(n)}}{\ud a}|_{a}\ud a\partial_{s}W(s,d(a))  \notag \\
&=&\int\limits_{a=0}^{t}\int\limits_{s=0}^{a}\mu ^{(n)}(s,d^{(n)}(a))\frac{%
\ud d^{(n)}}{\ud a}|_{a}\ud a\ud s+\int\limits_{a=0}^{t}\int\limits_{s=0}^{a}\Phi
^{(n)}(s,d^{(n)}(a))\partial_{s}W(s,d(a))\frac{\ud d^{(n)}}{\ud a}|_{a}\ud a  \notag \\
&=&\int\limits_{a=0}^{t}(x^{(n)}(a,d^{(n)}(a)))^{2}\frac{\ud d^{(n)}}{\ud a}|_{a}\ud a
\label{atC} \\
&=&\sum_{i=1}^{M_d^{(n)}}\sum_{k=1}^{N_{up}(i)}\int%
\limits_{s=J_{k,up}^{(n)}(i)}^{J_{next}(J_{k,up}^{(n)}(i))}(\tilde{x}_{i}^{(n)}(s))^{2}\frac{\ud d^{(n)}}{\ud s}|_{s}\ud s+\sum_{p=1}^{N_{down}(i)}\int\limits_{s=J_{p,down}^{(n)}(i)}^{J_{next}(J_{k,down}^{(n)}(i))}(\tilde{x}%
_{i}^{(n)}(s))^{2}\frac{\ud d^{(n)}}{\ud s}|_{s}\ud s.
\end{eqnarray}%
\bigskip

\noindent \textbf{Calculation of (A)}

We split $(A)$ in two parts:%
\begin{equation*}
(A)=(A1)+(A2)
\end{equation*}

where%
\begin{eqnarray*}
(A1) &=&\int\limits_{s=0}^{t}\int\limits_{u=0}^{d^{(n)}(s)}(-b(s)(\tilde{x}%
_{i}^{(n)}(s))^{2}-\frac{\sigma ^{2}(s)}{4})\ud u\ud s \\
&=&\int\limits_{s=0}^{t}\left( -b(s)Z_{mid}^{(n,M)}(s)-\frac{\sigma
(s)^{2}d^{(n)}(s)}{4}\right) \ud s.
\end{eqnarray*}

The calculation of (A2) is more complicated. It is easier to isolate the up
branches and the down branches of $d^{(n)}$:%
\begin{equation*}
(A2)=\sum_{i=1}^{M_d^{(n)}}\sum\limits_{k=1}^{N_{up}^{(n)}(i)}(A2_{k,up,i,t})+%
\sum\limits_{p=1}^{N_{down}^{(n)}(i)}(A2_{p,down,i,t})
\end{equation*}%
\bigskip

where we define:%
\begin{eqnarray*}
(A2_{k,up,i,t}) &=&-\int\limits_{s=J_{k,up}^{(n)}(i)}^{J_{next}^{(n)}(J_{k,up}^{(n)}(i))}\int\limits_{u=0}^{d^{(n)}(s)}\delta
(s-J_{k,up}^{(n)}(i))(\tilde{x}_{i}^{(n)}(s))^{2}1\{i-1<u\leq i\}\ud u\ud s
\label{A2kup} \\
&=&-\int\limits_{s=J_{k,up}^{(n)}(i)}^{J_{next}^{(n)}(J_{k,up}^{(n)}(i))}(\tilde{x}_{i}^{(n)}(J_{k,up}^{(n)}(i)))^{2}\frac{%
\ud d^{(n)}}{\ud s}|_{s}\ud s; \\
(A2_{p,down,i,t}) &=&-\int\limits_{s=J_{p,down}^{(n)}(i)}^{J_{next}^{(n)}(J_{p,down}^{(n)}(i))}\int%
\limits_{u=0}^{d^{(n)}(s)}\delta (t-J_{p,down}^{(n)}(i))(\tilde{x}%
_{i}^{(n,0,M)}(s))^{2}\frac{\ud d^{(n)}}{\ud s}|_{s}1\{i-1<u\leq i\}\ud u\ud s
\label{A2pdown} \\
&=&-\int\limits_{s=J_{k,up}^{(n)}(i)}^{J_{next}^{(n)}(J_{p,down}^{(n)}(i))}(\tilde{x}_{i}^{(n)}(J_{p,down}^{(n)}(i)))^{2}\frac{%
\ud d^{(n)}}{\ud s}|_{s}\ud s.
\end{eqnarray*}%
\bigskip

\noindent\textbf{Calculation of (B)}

By Levy's theorem there exists a collection of Brownian motions $B^{(n)}$
such that:%
\begin{eqnarray*}
&&(B)=\int\limits_{s=0}^{t}\int\limits_{u=0}^{d^{(n)}(s)}\sigma (s)\tilde{x}%
_{i}^{(n)}(s,\omega )\sum_{i=1}^{\infty }1\{i-1<u\leq i\}\ud u\ud W_{i}(s)
\label{B_up} \\
&&=\int\limits_{s=0}^{t}\sum\limits_{i=1}^{[d^{(n)}(s)]}\sigma (s)\tilde{x}%
_{i}^{(n)}(s)\ud W_{i}(s)+(d^{(n)}(s)-[d^{(n)}(s)])\sigma (s)\tilde{x}%
_{[d^{(n)}(s)]+1}^{(n)}(s)\ud W_{[d^{(n)}(s)]+1}(s)  \notag \\
&&=\int\limits_{s=0}^{t}\sigma (s)\sqrt{Z_{mid}^{(n)}(s)}\frac{%
\sum\limits_{i=1}^{[d^{(n)}(s)]}\tilde{x}%
_{i}^{(n)}(s)\ud W_{i}(s)+(d^{(n)}(s)-[d^{(n)}(s)])\tilde{x}%
_{[d^{(n)}(s)]+1}^{(n)}(s)\ud W_{[d^{(n)}(s)]+1}(s)}{\sqrt{\sum%
\limits_{i=1}^{[d^{(n)}(s)]}\tilde{x}_{i}^{(n)}(s)^{2}+\left(
(d^{(n)}(s)-[d^{(n)}(s)])\tilde{x}_{[d^{(n)}(s)]+1}^{(n)}(s)\right) ^{2}}}
\notag \\
&&=\int_{s=0}^{t}\sigma (s)\sqrt{Z_{mid}^{(n)}(s)}\ud B^{(n)}(s). \notag
\end{eqnarray*}%
Finally Lemma \ref{lemma2.2} follows from (\ref{ZABC}) and regrouping $(A)=(A1)+(A2)$, $(B)$ and $(C)$.
\bigskip

\noindent \textbf{Proof of Lemma \ref{lemma2.4}}

We observe that:%
\begin{eqnarray*}
&&J_{next}^{(n)}(J_{k,up}^{(n)}(i))-J_{k,up}^{(n)}(i)=o(\frac{1}{n}) \\
&&J_{next}^{(n)}(J_{p,down}^{(n)}(i))-J_{p,down}^{(n)}(i)=o(\frac{1}{n}).
\end{eqnarray*}

Clearly $\frac{\ud d^{(n)}}{\ud s}|_{s}=o(n)$. There is no jump of $\tilde{x}%
_{i,mid}^{(n,m,M)}$ in the interval $[J_{k,up}^{(n)}(i)\leq
s<J_{p(k),down}^{(n)}(i))$ if $s\in \lbrack t_{m}^{M},t_{m+1}^{M}]$. A
fortiori there is no jump of $\tilde{x}_{i,mid}^{(n,m,M)}$ in the interval $%
[J_{k,up}^{(n)}(i)\leq s<J_{next}^{(n)}(J_{k,up}^{(n)}(i)]$ if $s\in \lbrack
t_{m}^{M},t_{m+1}^{M}]$. Thus $E[|(\tilde{x}_{i,mid}^{(n,m,M)}(s))^{2}-(%
\tilde{x}_{i}^{(n)}(J_{k,up}^{(n)}(i)))^{2}|]$ is of order $1/M$ on $s\in $ $%
[L_{up}^{(m,M,k,i,n)},U_{up}^{(m,M,k,i,n)}]$, and:%
\begin{equation*}
E[\int\limits_{L_{up}^{(m,M,k,i,n)}}^{U_{up}^{(m,M,k,i,n)}}\left\vert (%
\tilde{x}_{i,mid}^{(n,m,M)}(s))^{2}-(\tilde{x}%
_{i,mid}^{(n)}(J_{k,up}^{(n,m,M)}(i)))^{2}\right\vert \ud s]=o(\frac{1}{M^{2}}).
\end{equation*}

Thus:%
\begin{eqnarray*}
E[\int\limits_{L_{up}^{(m,M,k,i,n)}}^{U_{up}^{(m,M,k,i,n)}}\left( (\tilde{x}%
_{i,mid}^{(n,m,M)}(s))^{2}\frac{\ud d^{(n)}}{\ud s}|_{s}-(\tilde{x}%
_{i,mid}^{(n,m,M)}(J_{k,up}^{(n)}(i)))^{2}\frac{\ud d^{(n)}}{\ud s}%
|_{J_{k,up}^{(n)}(i)}\right) \ud s] &\leq & \\
E[\int\limits_{L_{up}^{(m,M,k,i,n)}}^{U_{up}^{(m,M,k,i,n)}}\left\vert (%
\tilde{x}_{i,mid}^{(n,m,M)}(s))^{2}-(\tilde{x}%
_{i,mid}^{(n,m,M)}(J_{k,up}^{(n)}(i)))^{2}\right\vert \left\vert
\max\limits_{s\in \lbrack J_{k,up}^{(n)}(i),J^{next}(J_{k,up}^{(n)}(i))]}%
\frac{\ud d^{(n)}}{\ud s}|_{s}\right\vert \ud s] &=&o(\frac{n}{M^{2}})
\end{eqnarray*}

Since $Var[|(\tilde{x}_{i}^{(n,m,M)}(s))^{2}-(\tilde{x}%
_{i}^{(n,m,M)}(J_{k,up}^{(n)}(i)))^{2}|]$ is also of order $1/M$ on $s\in $ $%
[L_{up}^{(m,M,k,i,n)},U_{up}^{(m,M,k,i,n)}]$,

\begin{equation*}
Var[\int\limits_{L_{up}^{(m,M,k,i,n)}}^{U_{up}^{(m,M,k,i,n)}}\left( (\tilde{x%
}_{i,mid}^{(n,m,M)}(s))^{2}\frac{\ud d^{(n)}}{\ud s}-(\tilde{x}%
_{i,mid}^{(n,m,M)}(J_{k,up}^{(n)}(i)))^{2}\frac{\ud d^{(n)}}{\ud s}%
|_{J_{k,up}^{(n)}(i)}\right) \ud s]=o(\frac{n^{2}}{M^{2}})
\end{equation*}

The double sum $\sum_{i=1}^{M_d^{(n)}}\left(
\sum_{k=1}^{N_{up}(i)}+\sum_{k=1}^{N_{down}(i)}\right) $ contributes to a
number of jumps equal to $card\{$ $\mathcal{J}^{(n)}\}-1$, which is of order
$n$. Thus:%
\begin{equation}
E[\frac{R^{(n,m,M)}(t)}{n}]=o(\frac{n^{2}}{M^{2}})  \label{Ernm}
\end{equation}

The variance of $\frac{R^{(n)}(t)}{n}$ is a bit more difficult to calculate,
since \footnote{%
We focus only on the up branches, but the computation for the down branches
is similar.}(we focus only on the up branches, but the computation for the
down branches is similar):%
\begin{gather}
Var[\sum_{k=1}^{N_{up}^{(n)}(i)}\int%
\limits_{L_{up}^{(m,M,k,i,n)}}^{U_{up}^{(m,M,k,i,n)}}\left( (\tilde{x}%
_{i,mid}^{(n,m,M)}(s))^{2}\frac{\ud d^{(n)}}{\ud s}|_{s}-(\tilde{x}%
_{i,mid}^{(n)}(J_{k,up}^{(n,m,M)}(i)))^{2}\frac{\ud d^{(n)}}{\ud s}%
|_{J^{next}(J_{k,up}^{(n)}(i))}\right) \ud s]=  \notag \\
\sum_{k=1}^{N_{up}^{(n)}(i)}Var[\sum_{k=1}^{N_{up}^{(n)}(i)}\int%
\limits_{L_{up}^{(m,M,k,i,n)}}^{U_{up}^{(m,M,k,i,n)}}\left( (\tilde{x}%
_{i,mid}^{(n,m,M)}(s))^{2}\frac{\ud d^{(n)}}{\ud s}|_{s}-(\tilde{x}%
_{i,mid}^{(n,m,M)}(J_{k,up}^{(n)}(i)))^{2}\frac{\ud d^{(n)}}{\ud s}%
|_{J^{next}(J_{k,up}^{(n)}(i))}\right) \ud s]+  \notag \\
\sum_{\substack{ k,k^{\prime }=1  \\ k\neq k^{\prime }}}^{N_{up}^{(n)}(i)}%
\sum_{k=1}^{N_{up}^{(n)}(i)}Cov\left[ \int%
\limits_{L_{up}^{(m,M,k,i,n)}}^{U_{up}^{(m,M,k,i,n)}}\left( (\tilde{x}%
_{i,mid}^{(n,m,M)}(s))^{2}\frac{\ud d^{(n)}}{\ud s}|_{s}-(\tilde{x}%
_{i,mid}^{(n,m,M)}(J_{k,up}^{(n)}(i)))^{2}\frac{d^{{}}d^{(n)}}{\ud s}%
|_{J^{next}(J_{k,up}^{(n)}(i))}\right) \ud s,\right.  \label{Covar} \\
\left. \int\limits_{L_{up}^{(m,M,k^{\prime },i,n)}}^{U_{up}^{(m,M,k^{\prime
},i,n)}}\left( (\tilde{x}_{i,mid}^{(n,m,M)}(s))^{2}\frac{\ud d^{(n)}}{\ud s}|_{s}-(%
\tilde{x}_{i,mid}^{(n)}(J_{k^{\prime },up}^{(n,m,M)}(i)))^{2}\frac{\ud d^{(n)}}{%
\ud s}|_{J^{next}(J_{k^{\prime },up}^{(n)}(i))}\right) \ud s\right]  \notag
\end{gather}

However the number $N_{up}^{(n)}(i)$ of terms of $i$ is of order 1, thus the
number of terms in (\ref{Covar}) is still of order 1, and:%
\begin{equation}
Var[\frac{R^{(n,m,M)}(t)}{n}]=o(\frac{n^{2}}{M^{2}})  \label{Varrnm}
\end{equation}

Then by Markov's inequality: for all $\epsilon >0$ arbitrarily small,
\begin{equation*}
P\big(|\frac{R^{(n,m,M)}(t)}{M}|>\epsilon \big)\leq \frac{|\frac{%
R^{(n,m,M\,)}(t)}{M}|}{\epsilon }=\frac{\big(\frac{R^{(n,m,M)}(t)}{M}\big)}{%
\epsilon }\xrightarrow[n\rightarrow+\infty]{}0.
\end{equation*}%
As a result, $\frac{R^{(n,m,M)}(t)}{n}$ tends to $0$ in probability when $%
M\rightarrow \infty $.

\bigskip

\noindent \textbf{Proof of Lemma \ref{lemma2.5}} Fix any $\varepsilon ,\delta
>0$ . The proof is by induction on $n$. The case $n=1$ is trivial. Suppose
that it is true for $n-1$.

Let $A_{\varepsilon ,\delta }$ be the event such that $%
\{|r_{n,mid}^{(n,M)}(t_{m+1}^{M})-r_{n,mid}^{(n,M)}(t_{m+1}^{M}-\varepsilon
)|>\delta \}$. Set $B_{j,mid}^{(n)}=B_{j,mid}^{(n-1)}$ for $1\leq j\leq n-1$%
. Then:

\begin{equation*}
P(|%
\sum_{j=1}^{n-1}r_{j,mid}^{(n,M)}(t_{m+1}^{M})-r_{j,mid}^{(n,M)}(t_{m+1}^{M}-\varepsilon )|>\delta )=0.
\end{equation*}

Suppose $P(A_{\varepsilon ,\delta })>0$. We have then

\begin{equation*}
P(|%
\sum_{j=1}^{n}r_{j,mid}^{(n,M)}(t_{m+1}^{M})-r_{j,mid}^{(n,M)}(t_{m+1}^{M}-%
\varepsilon )|>\delta )>0.
\end{equation*}

By Levy's theorem, we arrive at the contradiction:

\begin{equation*}
P(|\sum_{j=1}^{n}Z_{mid}^{(n)}(t_{m+1}^{M})-Z_{mid}^{(n)}(t_{m+1}^{M}-%
\varepsilon )|>\delta )>0.
\end{equation*}

\noindent \textbf{Proof of Lemma \ref{lemma2.6}}

Define
\begin{equation*}
\tilde{R}^{(n,M)}(t)=\sum_{m=1}^M R^{(n,m,M)}(t)1\{t_{m}^{M}\leq t<t_{m+1}^{M}\}
\end{equation*}

\textbf{Step (i)}:\ tightness of $r_{1,mid}^{(n,M)}$

Let%
\begin{equation}
h_{1}^{(n,M)}(t)=\int_{s=0}^{t}-b(s)r_{1,mid}^{(n,M)}(s)-\frac{\sigma
^{2}(s)d(s)}{4}\ud s+\int_{s=0}^{t}\sigma (s)\sqrt{r_{1,mid}^{(n,M)}(s)}%
\ud B_{1,mid}^{(n,M)}(s).  \label{def_of_h1}
\end{equation}

Thus:%
\begin{equation*}
r_{1,mid}^{(n,M)}(t)=h_{1}^{(n,M)}(t)+\frac{\tilde{R}^{(n,M)}(t)}{n}.
\end{equation*}.
\bigskip

By Lemma \ref{lemma2.5}, the process $r_{1,mid}^{(n,M)}$ is continuous. Let the modulus
of continuity be:%
\begin{equation*}
w(r_{1,mid}^{(n,M)},\delta )=\sup_{\substack{ |s-t|<\delta  \\ 0\leq s,t\leq
T}}|r_{1,mid}^{(n,M)}(t)-r_{1,mid}^{(n,M)}(s)|.
\end{equation*}

We bound its expected value:%
\begin{equation*}
E[w(r_{1,mid}^{(n,M)},\delta )]\leq (A)+(B)+(C),
\end{equation*}

where:%
\begin{equation*}
(A)=E[\sup_{|s-t|<\delta }\int_{u=s}^{t}-b(u)r_{1,mid}^{(n,M)}(u)-\frac{%
\sigma ^{2}(u)d^{(n)}(u)}{4}\ud u]<K_{1}\delta
\end{equation*}

for some constant $K_{1}>0$. With the assumption that $\sigma $ is bounded:%
\begin{eqnarray*}
(B) &=&E[\sup_{|s-t|<\delta }\int_{u=s}^{t}\sigma (u)\sqrt{%
r_{1,mid}^{(n,M)}(u)}\ud B_{1,mid}^{(n,M)}(u)] \\
&\leq &E[\sup_{s\in \lbrack 0,T]}|\sigma (s)\sqrt{r_{1,mid}^{(n,M)}(s)}|\sup_{|s-t|<\delta }%
|B_{1,mid}^{(n,M)}(t)-B_{1,mid}^{(n,M)}(s)|] \\
&\leq &\big(E[\sup_{s\in \lbrack 0,T]}\sigma (s)\sqrt{r_{1,mid}^{(n,M)}(s)}%
]^{2}\big)^{1/2}\big(E[\sup_{|s-t|<\delta
}|B_{1,mid}^{(n,M)}(t)-B_{1,mid}^{(n,M)}(s)|]^{2}\big)^{1/2}~~%
\mbox{(H\"older
inequality)} \\
&\leq &K_{2}\delta ^{1/2},
\end{eqnarray*}

where $K_{2}>0$ is a constant. The last line follows because the expected
value of the maximum of Brownian motion $B$ is equal to :
\begin{equation*}
E[\sup_{s<t}|B(s)|]=\sqrt{\frac{2t}{\pi }}.
\end{equation*}

Finally%
\begin{equation*}
(C)=E\left[\sup_{\substack{ |s-t|<\delta  \\ 0\leq s,t\leq T}}|\frac{\tilde{R}%
^{(n,M)}(t)}{n}-\frac{\tilde{R}^{(n,M)}(s)}{n}|\right].  \label{C}
\end{equation*}

Since $r_{1,mid}^{(n,M)}$ is continuous at $t_{m}^{M}$ , let $%
t_{k-1}^{M}\leq s<t_{k}^{M}<..<\tau _{K}^{M}<t<t_{K+1}^{M}$. Then

\begin{eqnarray*}
|\frac{\tilde{R}^{(n,M)}(t)}{n}-\frac{\tilde{R}^{(n,M)}(s)}{n}| &<&|\frac{%
\tilde{R}^{(n,M)}(t_{k}^{M})}{n}-\frac{\tilde{R}^{(n,M)}(s)}{n}%
|+\sum_{i=k}^{K}|\frac{\tilde{R}^{(n,M)}(t_{i+1}^{M})}{n}-\frac{\tilde{R}%
^{(n,M)}(t_{i}^{M})}{n}|+ \\
&&|\frac{\tilde{R}^{(n,M)}(t)}{n}-\frac{\tilde{R}^{(n,M)}(\tau _{K}^{M})}{n}|
\\
&\leq &2\delta \max\limits_{k-1\leq i\leq K+1}|\frac{\tilde{R}%
^{(n,M)}(t_{i}^{M})}{n}|].
\end{eqnarray*}

Therefore

\begin{equation*}
\sup |\frac{\tilde{R}^{(n,M)}(t)}{n}-\frac{\tilde{R}^{(n,M)}(s)}{n}|\leq
2\delta \max\limits_{k-1\leq i\leq K+1}|\frac{\tilde{R}^{(n,M)}(t_{i}^{M})}{n%
}|.
\end{equation*}

Thus, by (\ref{Ernm}):

\begin{equation}
(C)=E[\sup_{\substack{ |s-t|<\delta  \\ 0\leq s,t\leq T}}|\frac{\tilde{R}%
^{(n,M)}(t)}{n}-\frac{\tilde{R}^{(n,M)}(s)}{n}|]=O(\frac{n^{2}\delta }{M^{2}}%
)
\end{equation}

Grouping the results, and passing to the sequence $M(n)=n^{4Q}$, we have:%
\begin{equation*}
E[w(r_{1,mid}^{(n,M(n))},\delta )]\leq K_{1}\delta +K_{2}\delta ^{1/2}
\end{equation*}

We first now demean the sequence $w(r_{1,mid}^{(n,M(n))},\delta )$. Clearly $%
E[w(r_{1,mid}^{(n,M(n))},\delta )]>0$. Then for any $\varepsilon >0$, using
Markov's inequality:
\begin{equation*}
P(w(r_{1,mid}^{(n,M(n))},\delta )\geq \varepsilon )\leq \frac{%
E[w(r_{1,mid}^{(n,M(n))},\delta )]}{\varepsilon }.
\end{equation*}

Clearly, for some constant $K_{4}>0:$%
\begin{equation*}
E[w(r_{1,mid}^{(n,M(n))},\delta )]\leq \sup_{|s-t|<\delta
}(E[(r_{1,mid}^{(n,M(n))}(t)-r_{j,mid}^{(n,M(n))}(s))^{2}])^{1/2}\leq
K_{5}\delta ^{1/2}.
\end{equation*}

Thus:%
\begin{equation*}
P(w(r_{1,mid}^{(n,M(n))},\delta )\geq \varepsilon )\leq \frac{K_{5}\delta
^{1/2}}{\varepsilon }.
\end{equation*}

We now invoke Theorem 8.2 in Billingsley (1968). Tightness occurs if, for each
positive $\varepsilon $ and $\eta $ there must exist of $\delta $ such that,
for all $n\geq n_{0}$%
\begin{equation}
P(w(r_{1,mid}^{(n,M(n))},\delta )\geq \varepsilon )\leq \eta.
\label{tightness}
\end{equation}

As a result we can select $\delta $ such that:%
\begin{equation*}
\frac{K_{5}\delta ^{1/2}}{\varepsilon }\leq \eta .
\end{equation*}

We thus select%
\begin{equation*}
\delta =\min (\frac{(\eta \varepsilon )^{2}}{K_{5}},1).
\end{equation*}

and, for all $n\geq n_{0}=\frac{T}{\delta }$, (\ref{tightness}) occurs.

\textbf{Step (ii)}: convergence of $r_{1,mid}^{(n,M(n))}$ to a solution of (\ref{CIR}).

By Lemma \ref{lemma2.4}, for any $\varepsilon >0$ and any $t$:%
\begin{equation*}
\lim_{n\rightarrow \infty }P(|r_{1,mid}^{(n,M(n))}(t)-h_{1}^{(n)}(t)|>\varepsilon
)=0
\end{equation*}

Thus,%
\begin{equation*}
r_{1,mid}^{(n,M(n))}(t)-h_{1}^{(n)}(t)\overset{d}{\rightarrow }0
\end{equation*}

Because of the Markov property, with $t_{1}<t_{2}:$%
\begin{equation*}
(r_{1,mid}^{(n,M(n))}(t_{1})-h_{1}^{(n)}(t_{1}),r_{1,mid}^{(n,M(n))}(t_{2})-r_{1,mid}^{(n,M(n))}(t_{1})-h_{1}^{(n)}(t_{2})-h_{1}^{(n)}(t_{1}))%
\overset{d}{\rightarrow }(0,0).
\end{equation*}

By Corollary 1 to Theorem 5.1 in Billingsley (1968):%
\begin{equation*}
(r_{1,mid}^{(n,M(n))}(t_{1})-h_{1}^{(n)}(t_{1}),r_{1,mid}^{(n,M(n))}(t_{2})-h_{1}^{(n)}(t_{2}))%
\overset{d}{\rightarrow }(0,0).
\end{equation*}

Repeating this argument, we can prove equality of the finite dimensional
distributions of $r_{1,mid}^{(n,M(n))}-h_{1}^{(n)}$ to zero. In other terms:%
\begin{equation*}
\lim_{n\rightarrow \infty }P(r_{1,mid}^{(n,M(n))}=h_{1}^{(n)})=1
\end{equation*}

We can extract a subsequence $\{n_{k(i)}\}$ so that $k\geq k(i)$ implies
that $P(|r_{1,mid}^{(n_{k},M(n_k))}-h_{1}^{(n_{k})}|\geq i^{-1})=2^{-i}$. By the
first Borel-Cantelli lemma there is a probability 1 that $%
|r_{1,mid}^{(n_{k})}-h_{1}^{(n_{k})}|\leq i^{-1}$ for all but finitely many $%
i$. Therefore:%
\begin{equation*}
\lim_{k\rightarrow +\infty }r_{1,mid}^{(n_{k},M(n_k))}=\lim_{k\rightarrow +\infty
}h_{1}^{(n_{k})}\text{ \ \ a.s. }P
\end{equation*}

In other terms:%
\begin{equation*}
P(r_{1,mid}^{(\infty,\infty )}=r_{1,mid}^{(\infty, \infty)}(0)-\int_{s=0}^{.}b(s)r_{1,mid}^{(\infty,\infty )}(s)-\frac{\sigma ^{2}d^{(n)}(s)%
}{4n}\ud s+\int_{s=0}^{.}\sigma (s)\sqrt{r_{1,mid}^{(\infty,\infty )}(s)}%
\ud B_{1,mid}^{(\infty )}(s))=1.
\end{equation*}

Thus $(r_{1,mid}^{(\infty,\infty )},B_{1,mid}^{(\infty )})$ is a weak solution to
the stochastic differential Equation (\ref{masterstoch2}). Maghsoodi (1996)
proved that under our conditions to (\ref{masterstoch2}) there is pathwise
uniqueness.

\bigskip

\textbf{Proof of Lemma \ref{lemma2.7}}

By differentiability of $d^{(n)}$, for any $t_{m}^{M}\leq t<t_{m+1}^{M}$:

\begin{equation*}
d^{(n)}(t)-d_{low}(t_{m}^{M},t_{m+1}^{M})=o(\frac{n}{M}).
\end{equation*}

By definition of $Z_{\chi ^{2}}^{(n,m,M)}$ and $\tilde{Z}_{mid}^{(n,M)}$

\begin{eqnarray}
E[Z_{\chi ^{2}}^{(n,m,M)}(t_{m+1}^{M})-Z_{mid}^{(n,M)}(t_{m+1}^{M})|\mathcal{%
F}_{t_{m}^{M}}] &=&o(\frac{n^{2Q}}{M^{2}};t_{m}^{M})  \label{EZ} \\
Var[Z_{\chi ^{2}}^{(n,m,M)}(t_{m+1}^{M})-Z_{mid}^{(n,M)}(t_{m+1}^{M})|%
\mathcal{F}_{t_{m}^{M}}] &=&o(\frac{n^{4Q}}{M^{2}};t_{m}^{M}).  \label{VarZ}
\end{eqnarray}

By differentiability of the characteristic function, the following Taylor series holds:

\begin{eqnarray*}
E[\exp (i\omega Z_{mid}^{(n,M)}(t))|\mathcal{F}_{t_{m}^{M}}] &=&E[\exp
(i\omega (Z_{\chi ^{2}}^{(n,m,M)}(t))|\mathcal{F}_{t_{m}^{M}}]+ \\
&&i\omega E[Z_{mid}^{(n,M)}(t)-Z_{\chi ^{2}}^{(n,m,M)}|\mathcal{F}%
_{t_{m}^{M}}] \\
&&-\frac{\omega ^{2}}{2}E[(Z_{mid}^{(n,M)}(t)-Z_{\chi ^{2}}^{(n,m,M)})^{2}|%
\mathcal{F}_{t_{m}^{M}}]+...
\end{eqnarray*}

Let $X$ be a $\mathcal{F}_{t_{m}^{M}}$-measurable random variable with
bounded mean. Suppose that
\begin{equation*}
X(\frac{K_{1}}{M^{2}}+i\frac{K_{2}}{M^{2}})\geq E[\exp (i\omega
Z_{mid}^{(n,M)}(t))|\mathcal{F}_{t_{m}^{M}}]-E[\exp (i\omega (Z_{\chi
^{2}}^{(n,m,M)}(t))|\mathcal{F}_{t_{m}^{M}}]\geq X(\frac{-K_{1}}{M^{2}}-i%
\frac{K_{2}}{M^{2}}),
\end{equation*}

then:%
\begin{gather}
|\omega E[Z_{mid}^{(n,M)}(t)-Z_{\chi ^{2}}^{(n,m,M)}|\mathcal{F}%
_{t_{m}^{M}}]-\frac{\omega ^{3}}{3!}E[(Z_{mid}^{(n,M)}(t)-Z_{\chi
^{2}}^{(n,m,M)})^{2}|\mathcal{F}_{t_{m}^{M}}]+..|\geq X\frac{K_{1}}{M^{2}}
\label{exp_o1} \\
|-\frac{\omega ^{2}}{2}E[(Z_{mid}^{(n,M)}(t)-Z_{\chi ^{2}}^{(n)})^{2}|%
\mathcal{F}_{t_{m}^{M}}]+\frac{\omega ^{4}}{3!}E[(Z_{mid}^{(n,M)}(t)-Z_{\chi
^{2}}^{(n,m,M)})^{2}|\mathcal{F}_{t_{m}^{M}}]+..|\geq X\frac{K_{2}}{M^{2}}
\label{exp_02}
\end{gather}

The left handside of (\ref{exp_o1}) and (\ref{exp_02}) can be majored by the
terms below, and we obtain
\begin{eqnarray*}
X\frac{K_{1}}{M^{2}} &\leq &\exp (\omega E[Z_{mid}^{(n,M)}(t)-Z_{\chi
^{2}}^{(n,m,M)}||\mathcal{F}_{t_{m}^{M}}])-1\leq X\frac{K_{3}}{M^{2}} \\
X\frac{K_{2}}{M^{2}} &\leq &\exp (\omega E[(Z_{mid}^{(n,M)}(t)-Z_{\chi
^{2}}^{(n,m,M)})^{2}|\mathcal{F}_{t_{m}^{M}}]))-1\leq X\frac{K_{4}}{M^{2}}
\end{eqnarray*}

which contradicts (\ref{EZ}) and (\ref{VarZ}). Thus

\begin{equation}
E[\exp (i\omega Z_{mid}^{(n,M)}(t))|\mathcal{F}_{t_{m}^{M}}]-E[\exp (i\omega
(Z_{\chi ^{2}}^{(n,m,M)}(t))|\mathcal{F}_{t_{m}^{M}}]=o(\frac{n^{4Q}}{M^{2}}%
;t_{m}^{M}). \label{approx_charac}
\end{equation}%
\bigskip \bigskip

From now on within this proof, we suppress the subscript $n$ from our
variables. For the moment we fix $M$, and take it out of the subscripts of
our variables. We let $\Delta =1/M$, and suppose that $K\Delta =t$. We find
it more convenient to shorten further the notation by:%
\begin{eqnarray*}
&&Z_{mid}^{(m,M)}(t_{m+1}^{M}) \rightarrow Z_{mid}((m+1)\Delta ) \\
&&d_{low}(m\Delta ,(m+1)\Delta ) \rightarrow \tilde{d}((m+1)\Delta ) \\
&&E[X|Z_{mid}^{(n)}(t_{m}^{M})] \rightarrow E_{m\Delta }[X]\\
&&c^{(n,m,M)}\rightarrow c(m\Delta ).
\end{eqnarray*}%
\bigskip

We denote the conditional non-centrality parameter and scale factor of $%
Z_{\chi ^{2}}((m+1)\Delta )$ given $Z_{mid}^{(n,M)}(t_{m}^{M})$ by:
\begin{eqnarray}
\lambda (m\Delta ) &=&c(m\Delta )Z_{mid}^{(n,M)}(m\Delta ))e^{-\tilde{b}%
(m\Delta )\Delta } \\
c(m\Delta ) &=&\frac{(1-e^{-\tilde{b}(m\Delta )\Delta })\tilde{\sigma}%
^{2}(m\Delta )}{4\tilde{b}(m\Delta )}.
\end{eqnarray}

where, so far $\tilde{b}$ and $\tilde{\sigma}$ are unknown functions. By the characteristic function of non-central chi-square distribution, the assumption that $M(n)=n^{4Q}$ and (\ref{approx_charac}), for $k=0,\ldots ,K-1$

\begin{equation}
E_{k\Delta }[\exp (i\omega Z_{mid}((k+1)\Delta ))]=\frac{\exp (i\frac{\omega
Z_{mid}(k\Delta )e^{-\tilde{b}(k\Delta )\Delta }}{1-2i\omega c(k\Delta )})}{%
(1-2i\omega c(k\Delta ))^{\tilde{d}((k+1)\Delta )/2}}+o(\Delta;k\Delta )
\label{z_chi}
\end{equation}%
As a way to find the general iteration formula, we study \footnote{%
We thank an anonymous referee for this judicious construction.} the case
where $k=2$. Defining:
\begin{equation*}
\omega _{2}=\frac{\omega e^{-\tilde{b}(2\Delta )\Delta }}{1-2i\omega
c(2\Delta )},
\end{equation*}

we calculate%
\begin{equation*}
E_{\Delta }[\exp (i\omega _{2}Z_{mid}(2\Delta ))]=\frac{\exp (i\frac{\omega
_{2}Z_{mid}(\Delta )e^{-\tilde{b}(\Delta )\Delta }}{1-2i\omega _{2}c(\Delta )%
})}{(1-2i\omega _{2}c(\Delta ))^{\tilde{d}(2\Delta )/2}}+o(\Delta;\Delta).
\end{equation*}

Hence,

\begin{eqnarray*}
E_{\Delta }[\exp (i\omega Z_{mid}(3\Delta ))] &=&\frac{\exp (i\frac{\omega
_{2}Z_{mid}(\Delta )e^{-\tilde{b}(\Delta )\Delta }}{1-2i\omega _{2}c(\Delta )%
})}{(1-2i\omega _{2}c(\Delta ))^{d_{low}(2\Delta )/2}(1-2i\omega c(2\Delta
))^{\tilde{d}(3\Delta )/2}}+o(\Delta ^{2};\Delta) \\
&=&\frac{\exp (i\frac{\omega Z_{mid}(\Delta )e^{(-\tilde{b}(\Delta )-\tilde{b%
}(2\Delta ))\Delta }}{1-2i\omega (c(2\Delta )+e^{-b(2\Delta )\Delta
}c(\Delta ))})}{(1-2i\omega (c(2\Delta )+e^{-\tilde{b}(2\Delta )\Delta
}c(\Delta )))^{\tilde{d}(2\Delta )/2}(1-2i\omega c(2\Delta ))^{(\tilde{d}%
(3\Delta )-\tilde{d}(2\Delta ))/2}}+o(\Delta;\Delta).
\end{eqnarray*}

Observe that, for a general case,
\begin{equation*}
E[\exp (i\omega Z_{mid}(K\Delta ))]=\frac{\exp \left( i\frac{\omega
Z_{mid}(0)e^{-\sum_{k=0}^{K-1}\tilde{b}(k\Delta )\Delta }}{1-2i\omega
(\sum_{k=0}^{K-1}e^{-\Delta \sum_{l=k+1}^{K-1}\tilde{b}(l\Delta )}c(k\Delta )%
}\right) }{(1-2i\omega \sum_{k=0}^{K-1}e^{-\Delta \sum_{l=k+1}^{K-1}\tilde{b}%
(l\Delta )}c(k\Delta ))^{\tilde{d}(\Delta )/2}A_{K}}+o(\Delta ),
\end{equation*}

where the denominator $A_{K}$ is given by:
\begin{equation*}
A_{K}=\prod_{j=0}^{K-1}\left( 1-2i\omega \sum_{k=j}^{K}e^{-\Delta
\sum_{l=k+1}^{K-1}\tilde{b}(l\Delta )}c(k\Delta )\right) ^{(\tilde{d}%
((j+1)\Delta )-\tilde{d}(j\Delta ))/2}.
\end{equation*}%
In order to show the limit of the above sequence, we take the logarithm and
get:
\begin{equation*}
\log A_{K}=\frac{1}{2}\sum_{j=0}^{K-1}(\tilde{d}((j+1)\Delta )-\tilde{d}%
(j\Delta ))\log \left( 1-2i\omega \sum_{k=j}^{K-1}e^{-\Delta
\sum_{l=k+1}^{K-1}\tilde{b}(l\Delta )}c(k\Delta )\right) .
\end{equation*}

We now return to placing subscripts $^{(M)}$ to our expressions. Let%
\begin{equation*}
\log A_{K/M}^{(M)}=\frac{1}{2}\sum_{j=0}^{K-1}h^{(M)}(\frac{j}{M})g^{(M)}(%
\frac{j}{M})
\end{equation*}%
\bigskip

where:%
\begin{eqnarray*}
h^{(M)}(\frac{j}{M}) &=&\tilde{d}(\frac{j+1}{M})-\tilde{d}(\frac{j}{M}) \\
&=&d_{low}^{{}}(\frac{j}{M},\frac{j+1}{M})-d_{low}^{{}}(\frac{j-1}{M},\frac{j%
}{M}) \\
g^{(M)}(\frac{j}{M}) &=&\log \left( 1-2i\omega \sum_{k=j}^{K-1}e^{-\frac{1}{M%
}\sum_{l=k+1}^{K-1}\tilde{b}(\frac{l}{M})}c^{(M)}(\frac{k}{M})\right)
\end{eqnarray*}

Let $\mathcal{R(}M\mathcal{)=\{}j|h(\frac{j}{M})\neq d^{\prime }(\frac{j}{M}%
)+o(\frac{1}{M}\}\}$. Since there is a finite number of minima of $d^{(n)}$, 
there is a sequence $M_K$ (take for instance a dyadic sequence) such that
for all $K>K_{0}$ the set $\mathcal{R}(M_K)$ is finite and constant. We
can thus split the calculations into

\begin{eqnarray}
&&\log A_{K}^{(M_K)}=\frac{1}{2}\sum_{j=0}^{K-1}h^{(M_K)}(\frac{j}{M_K}%
)g^{(M_K)}(\frac{j}{M_K})  \notag \\
&&=\sum_{j\in \{0,..,K\}\cap \overline{\mathcal{R(}M_K\mathcal{)}}}h^{(M_K)}(%
\frac{j}{M_K})g^{(M_K)}(\frac{k}{M_K})  \label{mon_deux} \\
&&+\sum_{j\in \mathcal{R(}M_K\mathcal{)}}h^{(M_K)}(\frac{j}{M_K})g^{(M_K)}(%
\frac{j}{M_K}).  \label{mon_trois}
\end{eqnarray}%
\bigskip

The sum in (\ref{mon_trois}) has a finite number of terms, while the sum in (%
\ref{mon_deux}) has a number of terms that tends to infinity when $%
K\rightarrow \infty $. 

Also, observe that, by using the mean value theorem,
\begin{equation*}
c(k\Delta )=\frac{(1-e^{-\tilde{b}(k\Delta )\Delta })\tilde{\sigma}%
^{2}(k\Delta )}{4\tilde{b}(k\Delta )}=\frac{\tilde{\sigma}^{2}(k\Delta
)\Delta }{4}+o(\Delta).
\end{equation*}
Thus we can write:

\begin{equation*}
\lim_{K\rightarrow \infty }\log A_{K}^{(M_K)}=\frac{1}{2}
\int_{0}^{T}d^{\prime }(s)g(s)\ud s,
\end{equation*}
with
$$
g(s)=\log (1-2i\omega
\int_{s}^{t}e^{-\int_{v}^{t}\tilde{b}(u)\,\mathrm{d}u}\tilde{\sigma}%
^{2}(v)/4\,\mathrm{d}v).
$$
Equivalently:%
\begin{equation*}
\lim_{\substack{ K\rightarrow \infty  \\ KM_{K}=t}}A_{K}^{(M_K)}=\exp (%
\frac{1}{2}\int_{0}^{t}d^{\prime }(s)\log (1-2i\omega
\int_{s}^{t}e^{-\int_{v}^{t}\tilde{b}(u)\,\mathrm{d}u}\tilde{\sigma}%
^{2}(v)/4\,\mathrm{d}v)\mathrm{d}s).
\end{equation*}
We also observe that, by Riemann sum,

$$\frac{\exp \left( i\frac{\omega
Z_{mid}(0)e^{-\sum_{k=0}^{K-1}\tilde{b}(k\Delta )\Delta }}{1-2i\omega
(\sum_{k=0}^{K-1}e^{-\Delta \sum_{l=k+1}^{K-1}\tilde{b}(l\Delta )}c(k\Delta )%
}\right) }{(1-2i\omega \sum_{k=0}^{K-1}e^{-\Delta \sum_{l=k+1}^{K-1}\tilde{b}%
(l\Delta )}c(k\Delta ))^{\tilde{d}(\Delta )/2}}\xrightarrow[K\rightarrow+\infty]{}
\frac{\exp (\frac{i\omega Z_{mid}(0)e^{-\int_0^t\tilde b(u)\ud u}}{1-2i\omega
\int_{0}^{t}e^{-\int_u^{t}\tilde b(u)\ud u}\tilde{\sigma} ^{2}(u)\ud u/4})}{(1-2i\omega \int_{0}^{t}e^{-\int_v^{t}\tilde b(u)\ud u}\tilde{\sigma}
^{2}(v)/4\ud v)^{d(0)/2}}.$$
Then Lemma \ref{lemma2.7} follows.
\bigskip

\textbf{Proof of Lemma \ref{lemma2.8}}

The proof is analogous to the proof of Lemma \ref{lemma2.7}, and we show only the first
inductive step. By (\ref{infinite_div}),

\begin{eqnarray*}
&&E[\exp (i\omega \tilde Z_{mid}^{(n,M)}(t_{m+1}^{M})|\mathcal{F}_{t_{m}^{M}}] \\
&&=E[\exp (i\omega \sum_{j=1}^{n}r_{j,mid}^{(n,M)}(t_{m+1}^{M})|\mathcal{F}%
_{t_{m}^{M}}] \\
&&= \exp (i\omega \sum_{j=1}^{n}r_{j,mid}^{(n,M)}(t_{m}^{M}))E[\exp (i\omega
\sum_{j=1}^{n}h_{j,mid}^{(n,M)}(t_{m+1}^{M})+R^{(n,m,M)})|\mathcal{F}%
_{t_{m}^{M}}] \\
&&= \exp (i\omega \tilde Z_{mid}^{(n,M)}(t_{m}^{M}))E[\exp (i\omega
\sum_{j=1}^{n}h_{j,mid}^{(n,M)}(t_{m+1}^{M})+R^{(n,m,M)})|\mathcal{F}%
_{t_{m}^{M}}].
\end{eqnarray*}

We use (\ref{Ernm}), (\ref{Varrnm}) as well as the arguments presented at the
beginning of Lemma \ref{lemma2.7} to obtain

\begin{eqnarray*}
&&E[\exp (i\omega \tilde Z_{mid}^{(n,M)}(t_{m+1}^{M})|\mathcal{F}_{t_{m}^{M}}] \\
&&= \exp (i\omega \tilde Z_{mid}^{(n,M)}(t_{m}^{M}))E[\exp (i\omega
\sum_{j=1}^{n}h_{j,mid}^{(n,M)}(t_{m+1}^{M})|\mathcal{F}_{t_{m}^{M}}]+o(%
\frac{n^{4Q}}{M^{2}};t_{m}^{M}) \\
&&= E\left[\exp \left( i\omega 
\sum_{j=1}^{n}\left(r_{j,mid}^{(n,M)}(t_{m}^{M})+h_{j,mid}^{(n,M)}(t_{m+1}^{M})%
\right) \right) |\mathcal{F}_{t_{m}^{M}}\right]+o(\frac{n^{4Q}}{M^{2}};t_{m}^{M}).
\end{eqnarray*}

However each copy of $%
r_{j,mid}^{(n,M)}(t_{m}^{M})+h_{j,mid}^{(n,M)}(t_{m+1}^{M})$ is
conditionally independent. Thus

$$
E[\exp (i\omega \tilde Z_{mid}^{(n,M)}(t_{m+1}^{M})|\mathcal{F}_{t_{m}^{M}}] = \left(E[\exp \left( i\omega \left(
r_{1,mid}^{(n,M)}(t_{m}^{M})+h_{1,mid}^{(n,M)}(t_{m+1}^{M})%
\right) \right) |\mathcal{F}_{t_{m}^{M}}]\right)^n+o(\frac{n^{2}}{M^{2}};t_{m}^{M}).
$$

Using the same method as before:

\begin{equation*}
E[\exp \left( i\omega \left(
r_{1,mid}^{(n,M)}(t_{m}^{M})+h_{1,mid}^{(n,M)}(t_{m+1}^{M})\right) \right) |%
\mathcal{F}_{t_{m}^{M}}]-E[\exp (i\omega r_{1,mid}^{(n,M)}(t_{m}^{M}))|%
\mathcal{F}_{t_{m}^{M}}]=o(\frac{n^{2Q}}{M^{2}};t_{m}^{M}).
\end{equation*}%
\bigskip

Thus:%
\begin{equation*}
E[\exp (i\omega \tilde Z_{mid}^{(n,M)}(t_{m+1}^{M})|\mathcal{F}_{t_{m}^{M}}]=
\Big(E[\exp (i\omega r_{1,mid}^{(n,M)}(t_{m+1}^{M}))|\mathcal{F}_{t_{m}^{M}}]\Big)^n+o(%
\frac{n^{4Q}}{M^{2}};t_{m}^{M}).
\end{equation*}
Then Lemma \ref{lemma2.8} holds by using the infinite divisibility of SNC chi-squared distribution and Lemma \ref{lemma2.7}. $\square$

\end{document}